\documentclass[12pt]{article}

\usepackage{hyperref}
\usepackage{amsthm,amsmath,amssymb}
\usepackage[latin1]{inputenc}
\usepackage[english]{babel}
\usepackage[a4]{crop}
\usepackage{graphics}
\usepackage{graphicx}
\usepackage{xypic}
\usepackage{multicol}

\newtheorem{theo}{Theorem}[section]

\newtheorem{prop}[theo]{Proposition}
\newtheorem{lemma}[theo]{Lemma}
\newtheorem{cor}[theo]{Corollary}
\theoremstyle{definition}

\newtheorem{ex}[theo]{Example}
\newtheorem{remark}[theo]{Remark}
\DeclareMathOperator{\tr}{tr}
\def\s{\sigma}
\def\l{\ell}
\def\Z{\mathbb{Z}}
\def\C{\mathbb{C}}
\def\j{\mathbb{\kappa}}
\def\root{\zeta}

\begin{document}

  \begin{center}
    \begin{Large}
      {\bf On the Varieties of Representations and Characters of 
      a Family of One-relator Groups. Extended Version} 
    \end{Large}

    \vspace{0.5cm}

    {\Large ($2^{\textrm{nd}}$ edition)}

    \vspace{0.5cm}

    {\large Jorge Mart\'{\i}n-Morales\footnote{{\em Email address:} {\tt jorge@unizar.es},
    {\em URL:} {\tt http://www.grupo.us.es/gmcedm/}} } \hspace{1cm}
    {\large Antonio M. Oller-Marc\'{e}n\footnote{{\em Email address:} {\tt oller@unizar.es},
    {\em URL:} {\tt http://www.unizar.es/matematicas/grupo\%20algebra/}}}\\[0.25cm]
    {\scriptsize\em Department of Mathematics-I.U.M.A. University of Zaragoza.}\\[-0.1cm]
    {\scriptsize\em C/ Pedro Cerbuna, 12 - 50009, Zaragoza (Spain).}
  \end{center}


  %
  %
  %
  %
  %
  %


  \begin{center}
  \begin{minipage}{14.5cm}
  \hrule
  \bigskip
  \noindent {\small {\bf Abstract}}

  \medskip
  \begin{footnotesize}
    \quad\ Let us consider the group $G = \langle x,y \mid
    x^m = y^n\rangle$ with $m$ and $n$ nonzero integers. In this paper, we study the variety of
    representations $R(G)$ and the character variety $X(G)$ in $SL(2,\C)$ of the group $G$, obtaining
    by elementary methods an explicit primary decomposition of the ideal corresponding to $X(G)$ in
    the coordinates $X=t_x$, $Y=t_y$ and $Z=t_{xy}$. As an easy consequence, a formula for computing
    the number of irreducible components of $X(G)$ as a function of $m$ and $n$ is given. We provide
    a combinatorial description of $X(G)$ and we prove that in most cases it is possible to recover
    $(m,n)$ from the combinatorial structure of $X(G)$. Finally we compute the number of irreducible
    components of $R(G)$ and study the behavior of the projection $t:R(G)\longrightarrow X(G)$.

  \end{footnotesize}
  \bigskip
  \hrule
  \end{minipage}
  \end{center}

  \section*{Introduction}

    Given a finitely generated group $G$, the set $R(G)$ of its representations over $SL(2,\C)$ can be
    endowed with the structure of an affine algebraic variety (see \cite{LM85}), the same holds for the
    set $X(G)$ of characters of representations over $SL(2,\C)$ (see \cite{CS83}). Since different
    presentations of a group $G$ give rise to isomorphic representation and character varieties; the
    study of geometric invariants of $R(G)$ and $X(G)$ like the dimension or the number of irreducible
    components is of interest in combinatorial group theory (see \cite{LI1,LI2,RUD} for instance). The
    varieties of representations and characters have also many applications in 3-dimensional geometry
    and topology (see \cite{HLM1,HLM2,TY} for instance). In \cite{HLM92} and \cite{HLM95} some aspects
    of the character variety of 2-bridge knots and links were studied, note that for a 2-bridge knot
    (or link) the fundamental group of its complement in $S^3$ admits a presentation with two generators
    and only one relation. In \cite{OLL} a geometrical description of the character variety of the torus
    knots $K_{m,2}$ was given.

    In \cite[Theor. 3.2.]{GM93} an explicit set of polynomials defining the character variety of a finitely
    presented group was given. Nevertheless this family of polynomials is not always satisfactory in order
    to give a geometrical description of the character variety. In this work, using elementary algebraic and
    arithmetic methods, we give an explicit primary decomposition of the ideal corresponding to the character
    variety of the group $G_{m,n}=\langle x,y\ |\ x^m=y^n\rangle$ with $m,n\neq0$, thus obtaining an easy
    geometrical description of it. This easier description allows us to compute geometrical invariants such
    as the number of irreducible components, to provide a combinatorial description of $X(G)$, and to study
    its relation with the variety of representations. The representation variety of $G_{m,n}$ was already
    studied in \cite{LI2} where its dimension and the number of irreducible components of maximal dimension
    were computed. We obtain here those results by a different approach and complete them by computing the
    number of irreducible components also for smaller dimension.

    Observe that if $\gcd(m,n)=1$ then $G_{m,n}$ is precisely the fundamental group of the exterior of the
    $(m,n)$-torus knot $K_{m,n}$, thus we have obtained the character variety for any torus knot. As an
    application we use our work to study the relation between the character variety of a torus knot $K_{m,n}$
    and that of its mirror image $K_{m,-n}$.

    The paper is organized as follows. In  Section \ref{XG}, we recall the construction of the character
    variety of a finitely presented group.
    In  Section \ref{sectionDefi}, we introduce some families of polynomials and give some technical results
    about them which are very useful in subsequent sections. Section \ref{mainResult} is devoted to give a
    complete description of the ideal associated with $X(G)$.
    In Section \ref{irrx}, we use this description to explicitly compute the number of irreducible
    components of $X(G)$ and in Section \ref{combXG} we study the combinatorial structure of $X(G)$.
    In Section \ref{inv} we study how much information about $m$ and $n$ is enclosed in this
    combinatorial structure.
    In Section \ref{irrR} we compute the number of irreducible components of the variety of
    representations $R(G)$ studying their relation with those of $X(G)$. Finally in Section \ref{com} we give
    some brief comments concerning the abelian component of $X(G)$, the relationship between the character variety
    of the exterior of a torus knot and that of its mirror image and an even simpler description of $X(G)$ in the
    case of the torus knots $K_{m,2}$.

  \section{Character Variety of Finitely Presented Groups}
  \label{XG}

    \bigskip
    Let $G$ be a group, a {\em representation} $\rho: G \longrightarrow
    SL(2,\C)$ is just a group homomorphism. We say that two representations $\rho$ and $\rho'$ are equivalent if
    there exists
    $P\in SL(2,\C)$ such that $\rho'(g)=P^{-1}\rho(g)P$ for every $g\in G$. A representation $\rho$ is
    {\em reducible} if the elements of $\rho(G)$ all share a
    common eigenvector, otherwise we say $\rho$ is {\em irreducible}. The following proposition presents some
    useful characterizations of reducibility.

    \begin{prop}\label{reducible} (see \cite[Lemma\ 1.2.1. and Prop.\ 1.5.5.]{CS83})
      \begin{enumerate}
      \item Let $\rho:G\longrightarrow SL(2,\C)$ be a representation. The following conditions are equivalent:
      \begin{enumerate}
        \item $\rho$ is reducible.
        \item $\rho (G)$ is, up to conjugation, a subgroup of upper triangular
          matrices.
        \item $\tr \rho(g) = 2$ for all $g$ in the commutator $G'=[G,G]$.
      \end{enumerate}
      \item If $G$ is generated by two elements $g$ and $h$, then $\rho:G\longrightarrow SL(2,\C)$ is
      reducible if and only if $\,\tr\rho ([g,h])=2$.
      \end{enumerate}
    \end{prop}

    Now, let us consider a finitely presented group
    $
      G = \langle x_1,\dots,x_k\ |\ r_1,\dots,r_s\rangle
    $
    and let $\rho:G\longrightarrow SL(2,\C)$ be a representation. It is clear that $\rho$ is completely
    determined by the $k$-tuple
    $(\rho(x_1),\dots,\rho(x_k))$ and thus we can identify
    $$
      R(G) = \{(\rho(x_1),\dots,\rho(x_k))\ |\ \rho\textrm{ is a representation of $G$}\}
      \subseteq\C^{4k}
    $$
    with the set of all representations of $G$ into $SL(2,\C)$,
    which is therefore (see \cite{CS83}) a well-defined affine algebraic set, up to
    canonical isomorphism.

    Recall that given a representation $\rho:G\longrightarrow SL(2,\C)$
    its character $\chi_{\rho}:G\longrightarrow \C$ is defined by
    $\chi_{\rho}(g)=\tr \rho(g)$. Note that two equivalent representations $\rho$ and $\rho'$ have
    the same character, and the converse is also true if $\rho$ or $\rho'$ is
    irreducible \cite[Prop. 1.5.2.]{CS83}. Now
    choose any $g\in G$ and define $t_g:R(G)\longrightarrow\mathbb{C}$ by
    $t_g(\rho)=\chi_{\rho}(g)$. It is easily seen that the ring $T$ generated by $\{t_g\ | \ g\in G\}$
    is a finitely generated ring (\cite[Prop.\ 1.4.1.]{CS83}) and, moreover, it
    can be shown using the well-known identities
    \begin{equation}\label{proptra}
    \begin{split}
    \tr A&=\tr A^{-1}\\
    \tr AB&=\tr BA\\
    \tr AB&=\tr A\tr B-\tr AB^{-1}
    \end{split}
    \end{equation}
    which hold in $SL(2,\mathbb{C})$ (see
    \cite[Cor. 4.1.2.]{GM93}) that $T$ is generated by the set:
    $$
      \{t_{x_i},t_{x_ix_j},t_{x_i x_j x_h}\ | \  1\leq i<j<h\leq k\}.
    $$

    Now choose $\gamma_1,\dots,\gamma_{\nu}\in G$ such that $T=\langle
    t_{\gamma_i}\ |\ 1\leq i\leq\nu\rangle$ and define the map
    $t:R(G)\longrightarrow\mathbb{C}^{\nu}$ by
    $t(\rho)=(t_{\gamma_1}(\rho),\dots,t_{\gamma_\nu}(\rho))$.
    Observe that $\nu \leq \frac{k(k^2+5)}{6}$.
    Put $X(G)=t(R(G))$, then $X(G)$ is an algebraic variety  which is well defined up
    to canonical isomorphism \cite[Cor.\ 1.4.5.]{CS83}
    and is called the \textit{character variety} of the group $G$ in $SL(2,\mathbb{C})$.
    Note that $X(G)$ can be identified with the set of all characters $\chi_\rho$
    of representations $\rho\in R(G)$.

    For every $0\leq j\leq k$ and for every $1\leq i\leq s$ we have that
    $p_{ij} = t_{r_ix_j}-t_{x_j}$ is a polynomial with rational
    coefficients in the variables $\{t_{x_{i_1}\dots x_{i_m}}\ |\
    m\leq3\}$, (see \cite[Cor. 4.1.2.]{GM93}). Then, we have the following explicit description of $X(G)$.
    \begin{theo}\label{theoremXG}
      (\cite[Theor.\ 3.2.]{GM93})
      $X(G)=\{\overline{x}\in X(F_k)\ |\
       p_{ij}(\overline{x})=0,\ \forall i,j\}$, where $F_k$ is the free group in $k$ generators.
    \end{theo}

    \begin{ex}\label{eje}
      \begin{enumerate}
      \item $X(F_1)=\C$,\quad $X(F_2)=\C^3$.
      \item Let $G=\langle x,y\ |\ xyx^{-1}y^{-1}\rangle$. It can be seen using the formulas
      given in (\ref{proptra}) (see \cite[Ex.\ 2]{HL05} for instance) that
      $X(G)=\{(X,Y,Z)\in\C^3\ |\ X^2+Y^2+Z^2-XYZ-4=0\}$. Observe that $G$ is the fundamental group of the
      two-dimensional torus. In what follows we will denote $D(X,Y,Z)=X^2+Y^2+Z^2-XYZ-4$, this polynomial will
      play a very important role in our paper since it satisfies $D(\tr A,\tr B,\tr AB) = \tr [A,B] - 2$ for
      all $A,B\in SL(2,\C)$.
      \end{enumerate}
    \end{ex}

    \begin{remark}
      Since the character variety of $X(F_1)$ is the whole field $\mathbb{C}$,
      the map $\tr:SL(2,\C) \longrightarrow \C$ given by the trace of a matrix is surjective. Therefore
      if two polynomials in one variable coincide on $\tr(SL(2,\C))$ then they are
      equal as polynomials.

      In the same way, since $X(F_2)= \C^3$, the map
      $t:SL(2,\C)\times SL(2,\C) \longrightarrow \C^3$
      given by $t(A,B) = (\tr A, \tr B, \tr AB)$ is surjective and thus if one wants
      to see that two polynomials in three variables are equal, it is enough to
      check that these polynomials coincide on $t(SL(2,\C)\times SL(2,\C))$.

      Note that, since $X(F_k)\neq\C^{\frac{k(k^2+5)}{6}}$ for all $k\geq3$, the previous considerations do
      not generalize for an arbitrary $k$. In the sequel we will make use of these observations without an
      explicit reference.
    \end{remark}

    We will end this section with a couple of technical results.
    \begin{lemma}\label{rax}
      Let $G$ and $H$ be groups and let $\varphi:R(G)\longrightarrow R(H)$ be a polynomial map.
      Then, there exists a unique polynomial map $\psi:X(G)\longrightarrow X(H)$ such that
      $t\circ\varphi=\psi\circ t$.
      $$
        \xymatrix{R(G) \ar[r]^{\varphi} \ar[d]_{t} \ar @{}[dr]|{\circlearrowleft} & R(H) \ar[d]^{t} \\
        X(G) \ar @{-->}[r]^{\psi} & X(H)}
      $$
    \end{lemma}

    \begin{cor}\label{inducedXG}
      Let $G$ and $H$ be groups and let $f:H \longrightarrow G$ be a group homomorphism.
      Then, there exists a unique polynomial map $\psi: X(G)\longrightarrow X(H)$ such that
      $\psi(t(\rho)) = t(\rho\circ f)$ for every $\rho\in R(H)$.
    \end{cor}

  \section{Some Families of Polynomials}
  \label{sectionDefi}

    In the forthcoming sections, we will need some particular families of polynomials whose
    definition and properties are given below.

    Given $c_0(T), c_1(T) \in \C[T]$ two polynomials, we define $\mathcal{F}_k^{(c_0,c_1)}(T)$
    as follows:
    $$
      \mathcal{F}_k^{(c_0,c_1)}: \begin{cases}
        \mathcal{F}_0^{(c_0,c_1)}(T) = c_0,\\
        \mathcal{F}_1^{(c_0,c_1)}(T) = c_1,\\
        \mathcal{F}_k^{(c_0,c_1)}(T) = T \mathcal{F}_{k-1}^{(c_0,c_1)}(T)
        - \mathcal{F}_{k-2}^{(c_0,c_1)}(T).
      \end{cases}
    $$
    Then we denote by $f_k$ and $h_k$ the polynomials $\mathcal{F}_k^{(2,T)}$
    and $\mathcal{F}_k^{(0,1)}$ respectively. Note that these families are closely related to the Chebyshev
    polynomials (see \cite{RI90}), in fact it can be shown that $f_k(2X)=2T_k(X)$
    and $h_k(2X)=U_{k-1}(X)$ for all $k$, where $T_k$ (resp. $U_k$) is the Chebyshev
    polynomial of the first (resp. second) kind.

    \smallskip
    Using (\ref{proptra}) it is possible to prove that
    $f_k$ is the only polynomial in one variable which satisfies
    $f_k(\tr A) = \tr A^k$ for all $A\in SL(2,\C)$. Analogously, given $a,b\in\Z$,
    it can be noted that there exists only one polynomial in three variables verifying
    $H(\tr A, \tr B, \tr AB) = \tr A^a B^{-b}$ for all $A, B\in SL(2,\C)$.
    We shall denote this polynomial by $F_{a,b}(X,Y,Z)\in \C[X,Y,Z]$. Now let us
    consider,
    $$
      s_k(T): \begin{cases}
        s_0(T) = 0,\\
        s_1(T) = s_2(T) = 1,\\
        s_3(T) = T+1,\\
        s_k(T) = T s_{k-2}(T) - s_{k-4}(T).
      \end{cases}
      \s_k(T): \begin{cases}
        \s_0(T) = 0,\\
        \s_1(T) = \s_2(T) = 1,\\
        \s_3(T) = T-1,\\
        \s_k(T) = T \s_{k-2}(T) - \s_{k-4}(T).
      \end{cases}
    $$
    Although we have defined these families of polynomials for $k\in\mathbb{N}$, they
    can clearly be extended to arbitrary $k\in \Z$.
    Finally, let $\j : \C^3 \to \C^3$ be the involution given by $\j(X,Y,Z) = (-X,-Y,Z)$.

    \begin{prop}\label{propibas}
    Let $a,b,i,j,k$ be integers, then:
      \begin{enumerate}
      \item $F_{a,b}(X,Y,Z) = F_{b,a}(Y,X,Z)$,\quad $F_{-a,-b}(X,Y,Z) = F_{a,b}(X,Y,Z)$.
      \item $F_{k,0}(X,Y,Z) = f_k(X)$,\quad $F_{0,k}(X,Y,Z) = f_k(Y)$.
      \item $s_{-k}(T) = -s_{k}(T)$,\ \ \ $\s_{-k}(T) = (-1)^{k-1}\s_{k}(T)$,\ \ \ $h_{-k}(T) = -h_k(T)$,\ \ \
        $f_{-k}(T) = f_k(T)$.
      \item $\j(f_k(X))=(-1)^kf_k(X)$,\quad $\j(s_k(X))=(-1)^{\left[\frac{k-1}{2}\right]}\s_k(X)$\quad and\quad
      $\j(F_{a,b}(X,Y,Z))=(-1)^{a-b}F_{a,b}(X,Y,Z)$.
      \item If $m$ is a positive integer, then
      \begin{footnotesize}
      \begin{multicols}{2}
      $
        s_m(T)= \begin{cases}
        1+\displaystyle{\sum_{i=1}^{\frac{m-1}{2}}f_i(T)} & \textrm{if $m$ is odd},\\
        \displaystyle{\sum_{i=1}^{\frac{m}{4}}f_{2i-1}(T)} & \textrm{if $m\equiv 0$},\\
        1+\displaystyle{\sum_{i=1}^{\frac{m-2}{4}}f_{2i}(T)} & \textrm{if $m\equiv 2$},\\
        \end{cases}
      $

      \columnbreak
      $
        \hspace{-0.5cm} (-1)^{\frac{m-1}{2}}\sigma_m(T) = 1+\displaystyle{\sum_{i=1}^{\frac{m-1}{2}}(-1)^i f_i(T)}\
        \textrm{ if $m$ is odd},
      $\\
      $
        h_m(T)= \begin{cases}
        1+\displaystyle{\sum_{i=1}^{\frac{m-1}{2}} f_{2i}(T)} & \textrm{if $m$ is odd},\\
        \displaystyle{\sum_{i=1}^{\frac{m}{2}} f_{2i-1}(T)} & \textrm{if $m$ is even},\\
        \end{cases}
      $
      \end{multicols}
      \end{footnotesize}
      where the congruences are taken modulo 4.

     \item $f_i(T)\cdot f_j(T) = f_{i+j}(T) + f_{i-j}(T)$.
     \item $h_{2k+1}(T) = s_{2k+1}(T) \s_{2k+1}(T)$,\quad
            $h_{2k}(T) = s_{2k}(T) f_{k}(T)$.
      \end{enumerate}
    \end{prop}

    \begin{proof}
    1) - 5) Follow by definition and/or inductive arguments.
    For 6) take $T=\tr A$ with $A\in SL(2,\C)$. Then $f_{i+j}(T) + f_{i-j}(T) =
    \tr A^i A^j + \tr A^i A^{-j} = \tr A^i \cdot \tr A^j = f_i(T)\cdot f_j(T)$.
    7) It is an easy consequence of 5) and 6) when $k$ is positive. The negative
    case follows immediately from the positive one together with 3) and the first part of 4).
    \end{proof}

  \section{Towards a Primary Decomposition}
  \label{mainResult}

    In what follows $m$ and $n$ will be assumed to be nonzero integers. Let $G_{m,n}$ be the group
    with presentation $G_{m,n} = \langle x, y \mid x^m = y^n \rangle$. Note that if $\gcd(m,n)=1$ this
    group is isomorphic to the fundamental group of the exterior of the $(m,n)$-torus knot. We are
    interested in computing $X(G_{m,n})$, its character variety in $SL(2,\C)$. Let $w=x^m y^{-n}$,
    then from Theorem \ref{theoremXG}, the ideal $J$ corresponding to $X(G_{m,n})$ can be generated by
    $t_{w}-2$, $t_{wx}-t_x$ and $t_{wy}-t_y$ in the ring of polynomials $\C[t_x,t_y,t_{xy}]$.
    In other words $J = \langle F_{m,n} (X,Y,Z) - 2, F_{m+1,n}(X,Y,Z) - X, F_{m,n-1}(X,Y,Z) - Y\rangle$,
    where $X=t_x$, $Y=t_y$ and $Z=t_{xy}$. Since $x^i y^{-k} = x^{j} y^{-l} \in G_{m,n}$ whenever $m=i-j$
    and $n=k-l$, all the polynomials $t_{x^i y^{-k}} - t_{x^j y^{-l}} = F_{i,k}(X,Y,Z) - F_{j,l}(X,Y,Z)$
    must belong to $J$. Therefore
    $
      J = \langle F_{i,k}(X,Y,Z) - F_{j,l}(X,Y,Z) \mid m = i-j, n=k-l\rangle
      \subset \C[X,Y,Z].
    $
    It is also possible to verify that $F_{i,k}-F_{j,l}\in J$ when $m=i+j$
    and $n=k+l$, and hence we can write, if necessary, $m=i+j$ and $n=k+l$ instead of
    $m=i-j$ and $n=k-l$ in the above expression of $J$.
    Now, associated with $(m,n)$, let us define some ideals in the
    ring $\C[X,Y,Z]$ which are closely related to $J$.
    $$
      \begin{array}{cc}
        I_1 = \langle s_m(X), s_n(Y) \rangle,\\
        \\
        I_3 = J + \langle D\rangle.
      \end{array}\quad
      I_2 =
      \begin{cases}
        \langle \s_m(X), \s_n(Y)\rangle \quad \text{if\quad $m,n$ are odd},\\
        \langle \s_m(X), f_{\frac{n}{2}}(Y)\rangle \quad \text{if\quad $m$ is odd and $n$ is even},\\
        \langle f_{\frac{m}{2}}(X), \s_n(Y)\rangle \quad \text{if\quad $m$ is even and $n$ is odd},\\
        \langle f_{\frac{m}{2}}(X), f_{\frac{n}{2}}(Y)\rangle \quad \text{if\quad $m,n$ are even},
      \end{cases}
    $$
    where $D(X,Y,Z)=X^2+Y^2+Z^2-XYZ-4$ is the polynomial defined in Example \ref{eje} 2).
    The aim of this section is to present the following theorem:

    \begin{theo}\label{mainTheorem}
      $V(J)=V(I_1\cap I_2\cap I_3) = V(I_1)\cup V(I_2) \cup V(I_3).$
    \end{theo}

    The proof of this theorem is rather technical and will be split into two parts.
    In the first part we will prove that $J$ is contained in $I_1 \cap I_2 \cap I_3$
    and in the second one we will prove the inclusion $I_1I_2I_3\subseteq J$. This suffices since we will have
    $I_1I_2I_3\subseteq J\subseteq I_1\cap I_2\cap I_3$, consequently $\sqrt{I_1I_2I_3}\subseteq
    \sqrt{J}\subseteq\sqrt{I_1\cap I_2\cap I_3}$ and, as we are working over an algebraically closed field,
    it follows that $\sqrt{I_1I_2I_3}=\sqrt{I_1\cap I_2\cap I_3}$ which completes
    the proof.

    \subsection{The First Inclusion}
    \label{mainResult1}

    In order to prove that $J$ is contained in $I_1 \cap I_2 \cap I_3$, we need
    some preliminary lemmas which give us some relations between $F_{a,b}$, $f_k$,
    $s_k$ and $\sigma_k$.



    \begin{lemma}\label{formula2}
      For all integers $a,b,i,j,k,l$, the following expressions hold:
      \begin{enumerate}
        \item $F_{a,k}(X,Y,Z) f_b(X) = F_{a+b,k}(X,Y,Z) + F_{a-b,k}(X,Y,Z)$.
        \item $F_{j,a}(X,Y,Z) f_b(Y) = F_{j,a+b}(X,Y,Z) + F_{j,a-b}(X,Y,Z)$.
        \item In particular, if $m=i-j$ is even then
          $F_{i,k}(X,Y,Z)+F_{j,k}(X,Y,Z)\in \langle f_{\frac{m}{2}}(X) \rangle$,
          and if $n=k-l$ is even then
          $F_{j,k}(X,Y,Z)+F_{j,l}(X,Y,Z)\in \langle f_{\frac{n}{2}}(Y) \rangle$.
      \end{enumerate}
    \end{lemma}

    \begin{proof}
      1) Let $A, B\in SL(2,\C)$ and put $(X,Y,Z) = (\tr A, \tr B, \tr AB)$. Then
      \setlength\multlinegap{0.4cm}
      \begin{multline*}
        F_{a+b,k}(X,Y,Z) + F_{a-b,k}(X,Y,Z) \ = \ \tr A^{a+b} B^{-k} + \tr A^{a-b} B^{-k} = \\
          = \tr A^b A^a B^{-k} + \tr A^{-b} A^a B^{-k} = \tr A^b\cdot \tr A^a B^{-k}
          \ = \ f_b(x) F_{a,k}(X,Y,Z).
      \end{multline*}

      2) Follows from 1), since $F_{a,b}(X,Y,Z) = F_{b,a}(Y,X,Z)$.

      3) Take $i=a+b$ and $j=a-b$ in 1). One obtains
      $$
        F_{i,k}(X,Y,Z) + F_{j,k}(X,Y,Z) = F_{\frac{i+j}{2},k}(X,Y,Z) f_\frac{m}{2}(X)
        \ \in\ \langle f_\frac{m}{2}(X) \rangle.
      $$
      In a similar way we can prove the second part of 3). In this case, take
      $k=a+b$ and $l=a-b$ in 2).
    \end{proof}

    \begin{lemma}\label{lemaGeneral}
      Let $\{p_k\}_{k\in\mathbb{Z}}$ be a family of polynomials satisfying the recursive equation
      $p_k = T p_{k-1} - p_{k-2}$
      for all $k\in\Z$ and let $\lambda,\mu :\Z\times\Z \to \Z$ be two maps verifying the following
      conditions:
      \begin{equation}\label{lambdaMu}
        \begin{cases}
          \lambda(i,j) = \lambda(i-1,j+1)\\
          \lambda(j+1,j) = j+1\\
          \lambda(j+2,j) = j+2
        \end{cases}\qquad
        \begin{cases}
          \mu(i,j) = \mu(i-1,j+1)\\
          \mu(j+1,j) = j\\
          \mu(j+2,j) = j
        \end{cases}\quad \forall i,j\in\Z.
      \end{equation}
      Then $p_i - p_j = s_{i-j}(T) \left(p_{\lambda(i,j)} - p_{\mu(i,j)}\right)$
      for all $i,j \in \Z$.
    \end{lemma}

    \begin{proof}
      First, suppose that $i-j\geq 0$ and let us prove the assertion by induction on
      $a=i-j\geq 0$. The claim is clear for $a=0,1,2$, since $s_0(T) = 0$,
      $s_1(T) = s_2(T) = 1$ and $\lambda, \mu$ satisfy the conditions in (\ref{lambdaMu})
      by hypothesis. For $a=3$ we note that $\lambda(j+3,j) = \lambda(j+2,j)=j+2$ and
      $\mu(j+3,j) = \mu(j+2,j+1) = j+1$, and hence
      \setlength\multlinegap{2.18cm}
      \begin{multline*}
        s_3(T) \left( p_{\lambda(j+3,j)} - p_{\beta(j+3,j)} \right) \ = \
          (T+1) ( p_{j+2} - p_{j+1} ) =\\
        (T p_{j+2} - p_{j+1}) - (T p_{j+1} - p_{j+2}) \ = \ p_{j+3} - p_j.
      \end{multline*}

      Let $a = i-j\geq 4$ and assume that the formula is true for all
      integers less than $a$ and greater than or equal to $0$. Then
      \setlength\multlinegap{0.5cm}
      \begin{multline*}
        s_{i-j}(T) \left( p_{\lambda(i,j)} - p_{\mu(i,j)} \right) \ = \
          \bigl( T s_{i-j-2}(T) - s_{i-j-4}(T) \bigl) \left( p_{\lambda(i,j)}
          - p_{\mu(i,j)} \right) =\\[0.15cm]
        = T s_{(i-1)-(j+1)}(T) \left( p_{\lambda(i-1,j+1)} - p_{\mu(i-1,j+1)}\right) -\hspace{5.8cm}\\
          \hspace{3.9cm} - s_{(i-2)-(j+2)}(T) \left( p_{\lambda(i-2,j+2)}
          - p_{\mu(i-2,j+2)} \right) \stackrel{(\text{IH})}{=}\\
        \stackrel{(\text{IH})}{=} T (p_{i-1}-p_{j+1}) - (p_{i-2} - p_{j+2}) = (Tp_{i-1} - p_{i-2})
          - (Tp_{j+1} - p_{j+2}) \ = \ p_i - p_j.
      \end{multline*}

      Now suppose that $a=i-j < 0$. Since $s_{-k}(T) = - s_k(T)$, to finish the proof,
      it is enough to see that $\lambda(i,j) = \lambda(j,i)$ and $\mu(i,j) = \mu(j,i)$.
      Assume, for instance, that $i > j$. There exists a positive integer $k$ such that
      $i = j+k$. Hence
      $$
        \lambda(i,j) = \lambda(j+k,j) = \lambda(j+k-1,j+1) = \cdots =
        \lambda(j,j+k) = \lambda(j,i).
      $$
      Analogously, $\mu(i,j) = \mu(j,i)$ and this finishes the proof.
    \end{proof}

    \begin{remark}
      In fact, $\lambda$ and $\mu$ are uniquely determined by conditions
      (\ref{lambdaMu}).
    \end{remark}

    \begin{cor}\label{formula3}
      Using the notation in Section \ref{sectionDefi}, we have:
      \begin{enumerate}
       \item $F_{i,k}(X,Y,Z) - F_{j,k}(X,Y,Z) = s_{i-j}(X)
         \left( F_{\left[\frac{i+j+2}{2}\right],k}(X,Y,Z) -
         F_{\left[\frac{i+j-1}{2}\right],k}(X,Y,Z) \right)$.
       \item $F_{i,k}(X,Y,Z) - F_{i,l}(X,Y,Z) = s_{k-l}(Y)
        \left( F_{i,\left[\frac{k+l+2}{2}\right]}(X,Y,Z) -
          F_{i,\left[\frac{k+l-1}{2}\right]}(X,Y,Z) \right)$.
      \end{enumerate}
      In particular, if $m=i-j$ then $F_{i,k}(X,Y,Z) - F_{j,k}(X,Y,Z)
      \in \langle s_m(X) \rangle$, and if $n=k-l$ then
      $F_{i,k}(X,Y,Z) - F_{i,l}(X,Y,Z) \in \langle s_{n}(Y) \rangle$.
      Moreover, if $m=i-j$ and $n=k-l$ are odd then $F_{i,k}(X,Y,Z) + F_{j,k}(X,Y,Z)$
      belongs to the ideal generated by $\s_m(X)$ while $F_{j,k}+F_{j,l}$ belongs to
      the ideal generated by $\s_n(Y)$.
    \end{cor}

    \begin{proof}
      1) Note that $\lambda(i,j) = [\frac{i+j+2}{2}]$ and $\mu(i,j) = [\frac{i+j-1}{2}]$
      satisfy the conditions (\ref{lambdaMu}) in Lemma \ref{lemaGeneral}, and we can
      take $p_i(X) = F_{i,k}(X,Y,Z)$ which verifies the corresponding recursive
      equation. Since $F_{a,b}(X,Y,Z) = F_{b,a}(Y,X,Z)$, $2)$ follows from $1)$.

      Now, assume that $m$ and $n$ are odd and let us apply to expression $1)$ the
      involution $\j:\C^3 \to \C^3$ given by $\j(X,Y,Z) = (-X,-Y,Z)$. Then the claim follows
      by Proposition \ref{propibas} $(4)$ and $(1)$.
    \end{proof}

    \begin{theo}\label{mainTheorem1}
      $J \subseteq I_1 \cap I_2 \cap I_3$.
    \end{theo}

    \begin{proof}
      Let $i,j,k,l$ be integers such that $m=i-j$ and $n=k-l$.
      Recall that $I_3$ is by definition the ideal $J+\langle D\rangle$.
      Therefore $J$ is contained in $I_3$. We divide the rest of
      the proof in two parts.

      \begin{enumerate}

        \item \underline{$J \subset I_1$}. The difference $F_{i,k}(X,Y,Z)-F_{j,l}(X,Y,Z)$
          can be written in the form
          $$
            F_{i,k}-F_{j,l} = (F_{i,k}-F_{j,k}) + (F_{j,k}-F_{j,l}),
          $$
          which belongs to $\langle s_m(X), s_n(Y) \rangle = I_1$ from Corollary \ref{formula3}.

        \item \underline{$J \subset I_2$}. Essentially, three cases can occur.
          \begin{itemize}
            \item $m$ and $n$ are odd. In this case, the involution $\j$ can be used.
              Since $J\subseteq I_1$, $\j(J) \subseteq \j(I_1)$. Now observe that
              from the fourth part of Proposition \ref{propibas}, $\j(J) = J$ and $\j(I_1) = I_2$.
            \item $m$ is odd and $n$ is even. Here, $F_{i,k}-F_{j,l}$ can be
              written in the form
              $$
                F_{i,k} - F_{j,l} = (F_{i,k} + F_{j,k}) - (F_{j,k} + F_{j,l})
              $$
              which belongs to the ideal $\langle \s_m(X), f_{\frac{n}{2}}(Y) \rangle$
              from Lemma \ref{formula2} and Corollary \ref{formula3}.
            \item $m$ and $n$ are even. We write the difference $F_{i,k}-F_{j,l}$
              as before. Now, from Lemma \ref{formula2}, it belongs to the ideal
              $\langle f_{\frac{m}{2}}(X), f_{\frac{n}{2}}(Y)\rangle$.
          \end{itemize}
      \end{enumerate}
    \end{proof}

    \subsection{The Second Inclusion}
    \label{mainResult2}

    We shall now show that $I_1 I_2 I_3 \subseteq J$.
    In order to do this, it is enough to study if
    $I_1 I_2 \langle D\rangle$ is contained in $J$, since $I_3 = J + \langle D\rangle$.
    If necessary, we will write $J_{m,n}$ instead of just $J$.
    Since $F_{a,b}(X,Y,Z) = F_{b,a}(Y,X,Z)$, if a polynomial
    $H(X,Y,Z)$ belongs to $J_{m,n}$ then $H(Y,X,Z)$ belongs to
    $J_{n,m}$. This easy remark allows us to simplify the
    proofs.

   Since $V(J)=X(G_{m,n})=\{(\tr A,\tr B, \tr AB)\ |\ A,B\in SL(2,\C),\ A^m=B^n\}$, we can assume
   $A^m=B^n$ when we work modulo $J_{m,n}$

    \begin{lemma}
      $h_m(X) D \in J_{m,n}$ for all $m, n\in \Z$. As a consequence,
      $h_n(Y) D \in J_{m,n}$ for all $m, n\in \Z$.
    \end{lemma}

    \begin{proof}
      Let $A, B$ be two matrices in $SL(2,\C)$ and put $(X,Y,Z)=(\tr A, \tr B, \tr AB)$.
      It can easily be proved by induction on $m$ that
      $h_m(X) D = \tr A^m B A^{-1} B^{-1} - f_{m-1}(X)$, since both members
      of the equality satisfy the same recursive equation.
      Now, recall that
      $A^m = B^n$ modulo $J_{m,n}$. Therefore $\tr A^m B A^{-1} B^{-1}
      \equiv \tr B^n A^{-1}$ and thus $h_m(X) D \equiv F_{1,n} - F_{m-1,0}\in J_{m,n}$.
    \end{proof}


    \begin{lemma}\label{lemmaSmD}
      Take $(X,Y,Z)\in \C^3$ and let $A, B\in SL(2,\C)$ be two matrices such that
      $(X,Y,Z) = (\tr A, \tr B, \tr AB)$. Then
      \begin{small}
      $$
        s_m(X) D =
        \begin{cases}
          \tr A^{\frac{m+1}{2}}BA^{-1}B^{-1}+\tr A^{\frac{m-1}{2}}BA^{-1}B^{-1}-\tr A^{\frac{m-3}{2}}
          -\tr A^{\frac{m-1}{2}} & \text{if $m$ is odd}, \\
          \tr A^{\frac{m}{2}}BA^{-1}B^{-1}-\tr A^{\frac{m-2}{2}} & \text{if $m$ is even}.
        \end{cases}
      $$
      \end{small}
    \end{lemma}

    \begin{proof}
      As before, the claim follows since both members of the expression
      verify the same recursive equation and it is obvious for $m=-1,0,1,2$.
    \end{proof}

    \begin{lemma}\label{demcomple}
      $s_m(X) f_{\frac{n}{2}}(Y) D \in J_{m,n}$ for all $m,n\in \Z$ with $n$ even.
      As a consequence, $f_{\frac{m}{2}}(X) s_n(Y) D \in J_{m,n}$ for all $m,n\in \Z$
      with $m$ even.
    \end{lemma}

    \begin{proof}
    We will make use of the previous lemma, so it is clear that we must work out
    separately the cases $m$ odd and $m$ even. We will only develop the odd case here.
    Let us recall that $f_{\frac{n}{2}}(Y)=\tr B^{\frac{n}{2}}$ and apply
    Lemma \ref{lemmaSmD} together with formulas (\ref{proptra}).
    \begin{equation}\label{Formula1}
      \begin{split}
        \quad f_{\frac{n}{2}}(Y)&s_m(X)D =\\[0.1cm]
        =&    \tr B^{\frac{n}{2}}\left(\tr A^{\frac{m+1}{2}}BA^{-1}B^{-1}
           +\tr A^{\frac{m-1}{2}}BA^{-1}B^{-1}-\tr A^{\frac{m-3}{2}}
           -\tr A^{\frac{m-1}{2}}\right) = \\[0.1cm]
        =&  \tr B^{\frac{n-2}{2}}A^{\frac{m+1}{2}}BA^{-1}
           +\tr A^{\frac{-m-1}{2}}B^{\frac{n-2}{2}}AB^{-1}
           +\tr B^{\frac{n-2}{2}}A^{\frac{m-1}{2}}BA^{-1}+\\
        &  +\tr A^{\frac{-m+1}{2}}B^{\frac{n+2}{2}}AB^{-1}-\tr B^{\frac{n}{2}}A^{\frac{m-3}{2}}
           -\tr B^{\frac{n}{2}}A^{\frac{-m+3}{2}}-\tr B^{\frac{n}{2}}A^{\frac{m-1}{2}}-\\
        &  -\tr B^{\frac{n}{2}}A^{\frac{-m+1}{2}}=\\[0.1cm]
        =&  \tr B^{\frac{n}{2}}\tr A^{\frac{m-1}{2}}B
           -\tr A^{\frac{-m-1}{2}}B^{\frac{n-2}{2}}AB^{-1}+
            \tr A^{\frac{-m-1}{2}}B^{\frac{n+2}{2}}AB^{-1}+\\
        &  +\tr B^{\frac{n}{2}}\tr A^{\frac{m-3}{2}}B
           -\tr A^{\frac{-m+1}{2}}B^{\frac{n-2}{2}}AB^{-1}
           +\tr A^{\frac{-m+1}{2}}B^{\frac{n+2}{2}}AB^{-1}-\\
        &  -\tr B^{\frac{n}{2}}A^{\frac{m-3}{2}}-\tr B^{\frac{n}{2}}A^{\frac{-m+3}{2}}
           -\tr B^{\frac{n}{2}}A^{\frac{m-1}{2}}-\tr B^{\frac{n}{2}}A^{\frac{-m+1}{2}}=\\[0.1cm]
       =&   \tr B^{\frac{n-4}{2}}A^{\frac{-m+1}{2}}
           -\tr A^{\frac{-m-1}{2}}B^{\frac{n-2}{2}}\tr AB^{-1}
           +\tr A^{\frac{-m-3}{2}}B^{\frac{n}{2}}+\\
        &  +\tr A^{\frac{-m-1}{2}}B^{\frac{n+2}{2}}\tr AB^{-1}
           -\tr A^{\frac{-m-3}{2}}B^{\frac{n+4}{2}}+\tr B^{\frac{n-4}{2}}A^{\frac{-m+3}{2}}-\\
        &  -\tr A^{\frac{-m+1}{2}}B^{\frac{n-2}{2}}\tr AB^{-1}
           +\tr A^{\frac{-m-1}{2}}B^{\frac{n}{2}}
           +\tr A^{\frac{-m+1}{2}}B^{\frac{n+2}{2}}\tr AB^{-1}-\\
        &  -\tr A^{\frac{-m-1}{2}}B^{\frac{n+4}{2}}-\tr B^{\frac{n}{2}}A^{\frac{-m+3}{2}}
           -\tr B^{\frac{n}{2}}A^{\frac{-m+1}{2}}=\\[0.1cm]
       =&  \left(\tr A^{\frac{-m-3}{2}}B^{\frac{n}{2}}
          -\tr A^{\frac{-m+3}{2}}B^{\frac{n}{2}}\right)
          +\left(\tr A^{\frac{-m-1}{2}}B^{\frac{n}{2}}
          -\tr A^{\frac{-m+1}{2}}B^{\frac{n}{2}}\right) +\\
        & +\left(\tr A^{\frac{-m+1}{2}}B^{\frac{n-4}{2}}
          -\tr A^{\frac{-m-1}{2}}B^{\frac{n+4}{2}}\right)
          +\left(\tr A^{\frac{-m+3}{2}}B^{\frac{n-4}{2}}
          -\tr A^{\frac{-m-3}{2}}B^{\frac{n+4}{2}}\right)+\\
        & +\tr AB^{-1}\left[\left(\tr A^{\frac{-m-1}{2}}B^{\frac{n+2}{2}}
          -\tr A^{\frac{-m+1}{2}}B^{\frac{n-2}{2}}\right)\right.+\\
        & +\left.\left(\tr A^{\frac{-m+1}{2}}B^{\frac{n+2}{2}}
          -\tr A^{\frac{-m-1}{2}}B^{\frac{n-2}{2}}\right)\right]=
    \end{split}
    \end{equation}

  \begin{equation*}
  \begin{split}
    =&  \left(F_{\frac{m+3}{2},\frac{n}{2}}-F_{\frac{m-3}{2}, \frac{n}{2}}\right)
       +\left(F_{\frac{m+1}{2},\frac{n}{2}}-F_{\frac{m-1}{2},\frac{n}{2}}\right) + \\
    &+  \left(F_{\frac{m-1}{2},\frac{n-4}{2}}-F_{\frac{m+1}{2}, \frac{n+4}{2}}\right)
       +\left(F_{\frac{m-3}{2},\frac{n-4}{2}}-F_{\frac{m+3}{2},\frac{n+4}{2}}\right) + \\
    &  + F_{1,1} \left[\left(F_{\frac{m+1}{2},\frac{n+2}{2}}-F_{\frac{m-1}{2},
    \frac{n-2}{2}}\right)+\left(F_{\frac{m-1}{2},\frac{n+2}{2}}-F_{\frac{m+1}{2},
    \frac{n-2}{2}}\right)\right].\\[0.25cm]
  \end{split}
  \end{equation*}

    Thus $s_m(X)f_{\frac{n}{2}}(Y)D\in J_{m,n}$ as claimed.
    \end{proof}

    \begin{lemma}
      $s_m(X) \s_n(Y) D \in J_{m,n}$ for all $m,n\in \Z$ with $n$ odd.
      As a consequence, $\s_m(X) s_n(Y) D \in J_{m,n}$ for all $m,n\in \Z$
      with $m$ odd.
    \end{lemma}

    \begin{proof}
      We will carry out the case when $m$ is odd and $n$ is positive (the case when $n$ is negative
      follows from the positive one and from Proposition \ref{propibas} (3) and (1)). First of all note that,
      $n$ being odd,
      $(-1)^{\frac{n-1}{2}}\s_n(Y) = 1+{\sum_{j=1}^{\frac{n-1}{2}}}(-1)^j \tr B^j$.
      Now applying again Lemma \ref{lemmaSmD} we can show, after some straightforward
      computations as in the previous lemma, that:
      \begin{equation}\label{Formula2}
        \begin{split}
          s_m(X)\s_n(Y)D = & \left(F_{\frac{m+3}{2},\frac{n-1}{2}}
          - F_{\frac{m-3}{2},\frac{n+1}{2}}\right) +
          \left(F_{\frac{m+1}{2},\frac{n-1}{2}}
          - F_{\frac{m-1}{2},\frac{n+1}{2}}\right) + \\
          & + \left(F_{\frac{m-1}{2},\frac{n+3}{2}}
          - F_{\frac{m+1}{2},\frac{n-3}{2}}\right) +
          \left(F_{\frac{m-3}{2},\frac{n+3}{2}}
          - F_{\frac{m+3}{2},\frac{n-3}{2}}\right) +\\
          & + F_{1,-1} \left[\left(F_{\frac{m+1}{2},\frac{n-1}{2}}-F_{\frac{m-1}{2},
          \frac{n+1}{2}}\right) + \left(F_{\frac{m-1}{2},\frac{n-1}{2}}-F_{\frac{m+1}{2},
          \frac{n+1}{2}}\right)\right].
        \end{split}
      \end{equation}
      and thus $s_m(X)\s_n(Y)D\in J_{m,n}$ as claimed.
    \end{proof}

    The last two lemmas can also be proved by double induction on $m$ and $n$, since both
    members of the corresponding expression satisfy the same recursive equation and
    hence the proofs are reduced to the base cases. However, we did not proceed in that
    direction because finding those identities would require to work out all the
    computations as in the proof of Lemma \ref{demcomple}.
    We refer to the appendix for complete details on how to find and show those
    complicated formulas.

    \begin{theo}\label{mainTheorem2}
      $I_1 I_2 I_3 \subseteq J$.
    \end{theo}

    \begin{proof}
      The above two lemmas imply that $I_1 I_2\langle D\rangle \subseteq J$. Now since
      $I_3 = J + \langle D\rangle$, we have that $I_1 I_2 I_3 = I_1 I_2 (J+\langle D\rangle) =
      I_1 I_2 J + I_1 I_2\langle D\rangle \subseteq J$.
    \end{proof}

    As we already remarked, Theorems \ref{mainTheorem1} and \ref{mainTheorem2} obviously imply Theorem
    \ref{mainTheorem}. Based on computational evidence, the authors conjecture that in fact
    $J = I_1 \cap I_2 \cap I_3$, but unfortunately no proof has been found.

  \section{Irreducible Components of $X(G_{m,n})$}
  \label{irrx}

    From Theorem \ref{mainTheorem}, the character variety $X(G_{m,n})$
    can be decomposed into the algebraic
    sets $V(I_1)$, $V(I_2)$ and $V(I_3)$. Therefore, in order to
    obtain an explicit description of $X(G_{m,n})$, we need to factorize
    the polynomials $s_k(T)$, $\s_k(T)$ and $f_k(T)$, and to find a nicer
    expression for $I_3 = J + \langle D\rangle$.

    \subsection{Finding the Factorization of $s_k(T)$, $\sigma_k(T)$ and $f_k(T)$}

    In this section $k$ will be a nonzero integer. Let us recall the recursive
    definition of the cyclotomic polynomials $\{c_{\l}(T)\}_{\l\in\mathbb{N}}$ given by
    \begin{equation}\label{ciclo}
    \prod_{\l|k}c_{\l}(T)=T^k-1,
    \end{equation}
    where $c_{\l}(T)$ is the minimal polynomial of the primitive $\l$-th roots of unity.
    In the same way we will denote
    by $r_{\l}(T)$ the minimal polynomial of the $\l$-th primitive roots of $-1$. This
    definition gives an expression similar to (\ref{ciclo}):
    \begin{equation}
    \prod_{\substack{\l|k\\ \frac{k}{\l}\ \textrm{odd}}}r_{\l}(T)=T^k+1.
    \end{equation}
    In fact it can be seen that if $\l$ is odd $r_{\l}(T)=c_{\l}(-T)$ and if $\l$ is even
    $r_{\l}(T)=c_{2\l}(T)$.

    Note that for every $3\leq \l\in\mathbb{N}$, there exists a polynomial (irreducible over $\Z$) $q_{\l}(T)$ such
    that: $$\displaystyle{c_{\l}(T)=T^{\frac{\varphi(\l)}{2}}q_{\l}\left(T+\frac{1}{T}\right)},$$
    where $\varphi$ is Euler's phi function and $\varphi(\l)=\textrm{deg}\ c_{\l}(T)$.
    For the sake of completeness, we will denote $q_1(T)=T-2$ and $q_2(T)=T+2$.
    \begin{lemma}\label{raizq}
    If $\l\geq 3$, $q_{\l}(T)$ has $\frac{\varphi(\l)}{2}$ distinct real roots and
    its set of zeroes is $Z[q_{\l}]=\{2\textrm{Re}\ z\ |\ z\
    \textrm{is a primitive $\l$-th root of unity}\}$ for all $\l\in\mathbb{N}$.
    Moreover if $\l_1\neq \l_2$, then $Z[q_{\l_1}]\cap Z[q_{\l_2}]=\emptyset$.
    \end{lemma}

    We can now obtain a factorization of the polynomial $f_k$ which will allow us to find its roots.

    \begin{prop}\label{fkt+t1}
    $\displaystyle{f_k(T)=\prod_{\substack{\l|k\\ \frac{k}{\l}\ \textrm{odd}}}q_{4\l}(T)}$.
    \end{prop}
    \begin{proof}
    \begin{equation*}
    \begin{split}f_k\left(T+\frac{1}{T}\right)&=T^k+\frac{1}{T^k}=\frac{T^{2k}+1}{T^k}=\frac{1}{T^k}\prod_{\substack{\l|2k\\ \frac{2k}{\l}\ odd}}r_{\l}(T)=\frac{1}{T^k}\prod_{\substack{\l|2k\\ \frac{2k}{\l}\ odd}}T^{\varphi(\l)}q_{2\l}\left(T+\frac{1}{T}\right)=\\ &=\prod_{\substack{\l|2k\\ \frac{2k}{\l}\ odd}}q_{2\l}\left(T+\frac{1}{T}\right)=\prod_{\substack{\l|k\\ \frac{k}{\l}\ odd}}q_{4\l}\left(T+\frac{1}{T}\right),
    \end{split}\end{equation*}
    where the last equality follows readily from the fact that $\{2\l\ \mid\ \l|2k,\
    \frac{2k}{\l}\ \textrm{odd}\}=\{4\l\ \mid\ \l|k,\ \frac{k}{\l}\ \textrm{odd}\}$ and
    from the identity $\displaystyle{\sum_{\l|2k,\, \frac{2k}{\l} \text{ odd}}\hspace{-0.35cm}\varphi(\l)=k}$.
    \end{proof}

    \begin{cor}
    The polynomial $f_k(T)$ has $k$ distinct real roots and its set of zeroes is $Z[f_k]=\{2\textrm{Re}\ z\ |\
    \textrm{z is a primitive $4\l$-th root of unity},\ \frac{k}{\l}\ \textrm{odd}\}=\{2\textrm{Re}\ z\ |\ z^{2k}=-1\}.$
    \end{cor}
    \begin{proof}
    The first equality is an immediate consequence of the previous proposition and Lemma \ref{raizq}.
    The second one can be verified by direct calculations.
    \end{proof}

    Now, we turn to the polynomials $s_k$ and $\s_k$. First we present a technical lemma that can easily be proved by induction.

    \begin{lemma}\label{skt+t1}
    $$s_k\left(T+\frac{1}{T}\right)=\begin{cases} \displaystyle{\frac{T^k-1}{T^{\frac{k-1}{2}}{(T-1)}}} & \textrm{if
    $k$ is odd.}\\ \displaystyle{\frac{T^k-1}{T^{\frac{k-2}{2}}{(T^2-1)}}} & \textrm{if $k$ is even.}\end{cases}$$
    \end{lemma}

    Using this result together with (\ref{ciclo}) and Lemma \ref{raizq}, and recalling Proposition \ref{propibas}, we have the following result.

    \begin{prop}\label{factors}
    $\displaystyle{s_k(T)=\prod_{1,2\neq \l|k}q_{\l}(T)},\quad
    \displaystyle{\s_k(T)=(-1)^{\left[\frac{k-1}{2}\right]}\prod_{1,2\neq \l|k}q_{\l}(-T)}$.\\
    Hence, $s_k(T)$ and $\s_k(T)$ have $\left[\frac{k-1}{2}\right]$ distinct real roots. In particular their sets of
    zeroes are $Z[s_k]=\{2\textrm{Re}\ z\ |\ z\neq \pm 1;\quad z^k=1\}=-Z[\s_k]$.
    \end{prop}

    \begin{remark}\label{signok}
    In so far $k$ was implicitly assumed to be positive. Observe that if we shift the sign of $k$ the sets of zeroes described above remain invariant. Namely
    $Z[f_k]=Z[f_{-k}]$, $Z[s_k]=Z[s_{-k}]$ and $Z[\s_k]=Z[\s_{-k}]$. Moreover, if $k$ is negative, then $f_k(T)$ has $|k|$ distinct real roots while $s_k(T)$ and $\s_k(T)$ have $\left[\frac{|k|-1}{2}\right]$.
    \end{remark}

    As a consequence of the previous discussion we will be able to count the number of irreducible components
    of $V(I_1)$ and $V(I_2)$ which, due to the form of the ideals $I_1$ and $I_2$, are disjoint straight lines.

    \subsection{Another Description for V($I_3$)}\label{descri}
    Now, we are interested in computing the number of irreducible
    components of the affine algebraic variety associated with $I_3$.
    For this aim, our definition of $I_3 = J+\langle D\rangle$ doesn't seem to be very
    useful.
    That is why we need another expression for $\mathcal{C}:=V(I_3)$. From now on we assume that
    $d=\gcd (m,n)$ and we will write $m=m'd$ and $n=n'd$.

    \begin{lemma}\label{anotherVI3}
    $V(I_3) = \{ (u+u^{-1},v+v^{-1},uv+(uv)^{-1}) \mid u,v\in\C^{*},\ u^m = v^n \}$.
    \end{lemma}

    \begin{proof}
    Let $t:R(G_{m,n})\longrightarrow X(G_{m,n})$ be the polynomial map defined in Section \ref{XG} and
    let us recall that by definition $V(J) = X(G_{m,n}) = t(R(G_{m,n}))$.
    Since $D(t(\rho)) = \tr\rho[x,y]-2$ when $\rho\in R(G_{m,n})$
    and $V(I_3) = V(J+\langle D\rangle) = V(J)\cap V(D)$, it is clear that
    $
      V(I_3) = \{ t(\rho) \mid \rho\in R(G_{m,n}),\ \tr \rho[x,y]=2 \}.
    $
    Now the claim follows from Proposition \ref{reducible}(2) and the identity
    $t(Red)=t(Diag)$, where $Red$ is the set of all reducible representations in $R(G_{m,n})$ and
    $Diag$ the set of all diagonal ones, (see the proof of Corollary 1.4.5 in \cite{CS83}).
    \end{proof}

    Using the previous lemma, since $u^m-v^n = \prod_{i=0}^{d-1} (u^{m'}-\root^i v^{n'})$ where
    $\{1,\ldots,\root^{d-1}\}$ is the set of all $d$-th roots of unity, the variety
    associated with $I_3$ can be decomposed as the union $\bigcup_{i=0}^{d-1} C_{\root^i}$ with
    $C_{\root^i} = \{ (u+u^{-1},v+v^{-1},uv+(uv)^{-1}) \mid u,v\in\C^{*},\ u^{m'} = \root^i v^{n'} \}$. Once
    we remove the redundant components, the above union provides a decomposition of $\mathcal{C}=V(I_3)$ into irreducible
    components. This fact is clarified in the lemma below.

    \begin{lemma}\label{Czetai}
    The following properties hold:
    \begin{enumerate}
      \item $C_{\root^i}$ is an irreducible algebraic set for all $i$.
      \item $C_{\root^i} = C_{\root^j}$ if and only if $\root^i = \root^{\pm j}$.
      Moreover, if $C_{\root^i}\neq C_{\root^j}$ then they are disjoint.
    \end{enumerate}
    \end{lemma}

    \begin{proof}
    1) Let us first show that the set $C_{\root^i}$ is algebraic.
    Take $Diag \subset R(G_{m,n})$ the algebraic set of all diagonal representations and $Diag^{\root^i}$
    the algebraic subset
    \begin{small}
    $$
      Diag^{\root^i} = \left\{ \left[\left(\begin{array}{cc} u & 0\\ 0 & u^{-1}\end{array}\right),
      \left(\begin{array}{cc} v & 0\\ 0 & v^{-1}\end{array}\right)\right]\in SL(2,\C)^2 \mid
      u^{m'} = \root^i v^{n'} \right\}.
    $$
    \end{small}
    \hspace{-0.275cm} Since the map $\C^{*}\longrightarrow \C$ given by $a\mapsto a+a^{-1}$ is proper with the usual topology over $\C$, then so is $t|_{Diag}$ and, consequently, also $t|_{Diag^{\root^i}}$ and it follows (see \cite{Serre55} and \cite{Hartshorne77}) that $t_{|Diag^{\root^i}}$ is closed with the Zariski topology.
    Now observe that $C_{\root^i} = t(Diag^{\root^i})$ and $Diag^{\root^i}$ is Zariski-closed,
    (see also the proof of Corollary 1.4.5 in \cite{CS83}).
    Finally, the irreducibility of $C_{\root^i}$ follows from the irreducibility of $Diag^{\root^i}$
    which is birationally equivalent to the plane curve given by the equation $u^{m'} = \root^i v^{n'}$.

    2) Note that the equation $a+a^{-1} = b + b^{-1}$ with $a,b\in \C^{*}$ has the only solutions
    $a=b^{\pm 1}$. Now, the assertion is a clear consequence of this fact.
    \end{proof}

    Using a parametrization of the plane curve
    given by the equation $u^{m'} = \root^{i} v^{n'}$ we can also parametrize $C_{\root^i}$. In this way it is easy to compute the self-intersection and regular points of $C_{\root^i}$. Since the rational morphism $\C^{*} \longrightarrow \C: t\mapsto t+1/t$
    is a $2:1$ branched covering, it is also easy to prove
    that $C_{\root^i}$ is smooth.
    Finally, all of its components being
    smooth and disjoint, we conclude that $\mathcal{C}$ is smooth.

    \subsection{Counting the Irreducible Components of $X(G_{m,n})$}
    Finally we present the main result of this section, which allows us to explicitly
    count the number of irreducible components of $X(G_{m,n})$.

    \begin{theo}\label{compx}
    The number of irreducible components of $X(G_{m,n})$ is:
    $$\begin{cases}
    \displaystyle{\frac{(|m|-1)(|n|-1)}{2}+\frac{d+1}{2}} & \textrm{if $d$ is odd},\\
    \ &\ \\
    \displaystyle{\frac{(|m|-1)(|n|-1)+1}{2}+\frac{d+2}{2}} & \textrm{if $d$ is even},
    \end{cases}$$
    where the first summand corresponds to the number of straight lines in $V(I_1)\cup V(I_2)$
    and the second one corresponds to the number of irreducible components of $V(I_3)$.
    \end{theo}

    Note that if $d=1$ the genus of the torus knot $K_{m,n}$ is precisely the number of
    straight lines in its character variety.

    \section{Combinatorial Description of $X(G_{m,n})$}\label{combXG}
    In this section we are interested in studying the combinatorial structure
    of $X(G)$ as well as in describing the information that can be obtained from it.
    Recall that the character varieties of $G$ associated with $(m,n)$, $(m,-n)$
    and $(n,m)$ are all isomorphic, so from now on we will assume
    that $m\geq n$ are positive integers.

    We know that the straight lines of $\mathcal{L}:= V(I_1)\cap V(I_2)$ are disjoint
    and so are the irreducible components of $\mathcal{C}=V(I_3)$. Therefore in order to study the combinatorial structure of
    $X(G)=\mathcal{L}\cup\mathcal{C}$ it is enough to study how the straight lines and
    the irreducible components of $\mathcal{C}$ intersect. The following result will be useful in the sequel.

    \begin{lemma}
    $\# L\cap \mathcal{C} = 2$ for all $L\in\mathcal{L}$.
    \end{lemma}

    \begin{proof}
    Consider $L:\{x=a, y=b\}$; i.e., $a$ is a root of one of the polynomials $s_m$, $\sigma_m$ or $f_{\frac{m}{2}}$
    and $b$ is a root of one of the polynomials $s_n$, $\sigma_n$ or $f_{\frac{n}{2}}$ depending
    on the parity of $m$ and $n$. Now take $\lambda$ and $\mu$ such that $a=\lambda+\lambda^{-1}$ and
    $b=\mu+\mu^{-1}$. It is clear that $L\cap \mathcal{C}$ is contained in
    $\{ (a,b,\lambda \mu + \lambda^{-1} \mu^{-1}), (a,b,\lambda \mu^{-1} + \lambda^{-1} \mu) \}
    =\{P_1, P_2\}$. We need to prove that these two points are different and they indeed
    belong to $\mathcal{C}$. If $P_1 = P_2$ then $\lambda \mu + \lambda^{-1} \mu^{-1} =
    \lambda \mu^{-1} + \lambda^{-1} \mu$ and thus $\lambda = \pm 1$ or $\mu = \pm 1$.
    This contradicts Lemma \ref{skt+t1} and Proposition \ref{fkt+t1}.
    Now, note that $P_1, P_2\in \mathcal{C}$ if and only if $\lambda^m = \mu^n$ and
    $\lambda^m = \mu^{-n}$. The claim follows again from Lemma \ref{skt+t1} and
    Proposition \ref{fkt+t1}.
    \end{proof}

    \begin{remark}
    In fact it can be shown that if a line of the form $L:\{x=a, y=b\}$
    and $\mathcal{C}$ intersect exactly in two different points, then $L$ must be in $\mathcal{L}$.
    \end{remark}

    In the light of the previous lemma and given $L\in\mathcal{L}$ we see that $L$ can only intersect $\mathcal{C}$ in one or two of its components. It is clear
    that a point $(a,b,c)$ with $a=\lambda+\lambda^{-1}$ and $b=\mu+\mu^{-1}$ belongs to $C_{\root^i}$
    if and only if $\lambda^{m'} = \root^{\pm i} \mu^{\pm n'}$. This means that $L:\{x=a, y=b\}$
    intersects $\mathcal{C}$ only in the components $C_{\lambda^{m'}\mu^{n'}}$ and $C_{\lambda^{m'}\mu^{-n'}}$ (which can be equal or not). After this discussion and Lemma \ref{Czetai},
    the following result is satisfied.

    \begin{prop}\label{corte}
    Let $L:\{x=a, y=b\}$ be a straight line of $\mathcal{L}$ and let $\lambda, \mu$ be two numbers such
    that $a=\lambda+\lambda^{-1}$ and $b=\mu+\mu^{-1}$. Then $C_{\lambda^{m'}\mu^{n'}}$
    and $C_{\lambda^{m'}\mu^{-n'}}$ are the only components of $\mathcal{C}$ which intersect $L$.
    In particular $L$ intersects only one component of $\mathcal{C}$ at two different points if and only if
    $\lambda^{m'} = \pm 1$ or $\mu^{n'} = \pm 1$.
    \end{prop}

    We have just seen that a line $L\in\mathcal{L}$ intersects $\mathcal{C}$ in two different points which can lie either in the same irreducible component of $\mathcal{C}$ or in different ones. Now we will compute the number of straight lines in each situation.

  \subsection{Straight Lines Intersecting Only One Component}

We will start by computing the number of lines of $\mathcal{L}$ intersecting $\mathcal{C}$ in only one component. Of course this component can be $C_1$, $C_{-1}$ (only if $d$ is even) and $C_{\root^i}$ (with $i\neq0,\frac{d}{2}$).

    \begin{prop}
    The number of straight lines intersecting $C_1$ twice is $\frac{(m'-1)(n'-1)}{2}$.
    Moreover, the previous number can be decomposed as follows.
    $$
      \frac{(m'-1)(n'-1)}{2} = \left\{
      \begin{array}{ccccl}
        \frac{(m'-1)(n'-1)}{4} & + & \frac{(m'-1)(n'-1)}{4} &
        \ \text{if}\ & \text{$d, m', n'$ odd,}\\[0.15cm]
        \frac{(m'-2)(n'-1)}{4} & + & \frac{m'(n'-1)}{4} &
        \ \text{if}\ & \text{$d, n'$ odd, $m'$ even,}\\[0.15cm]
        \frac{(m'-1)(n'-1)}{2} & + & 0 &
        \ \text{if}\ & \text{$d$ even,}
      \end{array}\right.
    $$
    where the first (resp. second) summand corresponds to the number of straight lines in $V(I_1)$
    (resp. $V(I_2)$).
    \end{prop}
    \begin{proof}
    We will assume that $d$, $m'$ and $n'$ are odd. The other three possible cases follow from similar considerations (observe that $m'$ and $n'$ cannot be both even). Recall that $\mathcal{L}=V(I_1)\cup V(I_2)$, so we will count the lines from $V(I_1)$ and $V(I_2)$ separately.

    Let $L:\{x=a,y=b\}$ be a straight line from $V(I_1)$, then $a=\lambda+\lambda^{-1}$ and $b=\mu+\mu^{-1}$ with $\lambda^m=\mu^n=1$ and $a,b\neq\pm 2$. From Proposition \ref{corte} it follows that $\lambda^{m'}=\mu^{n'}=\pm 1$ and since $d$ is odd it must be $\lambda^{m'}=\mu^{n'}=1$. Now, the previous system of complex equations has $(m'-1)(n'-1)$ solutions in $(\lambda,\mu)$ (recall that $\lambda,\mu\neq\pm 1$) but we are interested in the possible values of $(a,b)$. Clearly $(\lambda^{\pm},\mu^{\pm})$ provide the same value $(a,b)$ thus, the number of lines from $V(I_1)$ intersecting $C_1$ twice is $\frac{(m'-1)(n'-1)}{4}$.

    On the other hand if $L:\{x=a,y=b\}$ is a straight line in $V(I_2)$, then $a=\lambda+\lambda^{-1}$ and $b=\mu+\mu^{-1}$ with $\lambda^m=\mu^n=-1$. Analogously to the previous paragraph, we conclude that the number of straight lines in $V(I_2)$ intersecting $C_1$ twice is again $\frac{(m'-1)(n'-1)}{4}$ and the proposition follows.
    \end{proof}

    \begin{prop}
    If $d$ is even, the number of straight lines intersecting $C_{-1}$ twice is $\frac{(m'-1)(n'-1)}{2}$.
    Moreover, all such straight lines are in $V(I_1)$.
    \end{prop}

    \begin{proof}
    We will only work out the case when $m'$ is even and $n'$ is odd, the other case being analogous. As in the previous proposition we will compute the lines coming from $V(I_1)$ and $V(I_2)$ separately.

    Let $L:\{x=a,y=b\}$ be a straight line in $V(I_1)$, then $a=\lambda+\lambda^{-1}$ and $b=\mu+\mu^{-1}$ with $\lambda^m=\mu^n=1$ and $a,b\neq\pm 2$. From Proposition \ref{corte} it follows that $\lambda^{m'}=-\mu^{n'}=\pm 1$. Now, the system of complex equations $\lambda^{m'}=1=-\mu^{n'}$ has $(m'-2)(n'-1)$ solutions in $(\lambda,\mu)$ while the system $\lambda^{m'}=-\mu^{n'}=-1$ has $m'(n'-1)$. Again, since we are interested in the possible values of $(a,b)$, it follows that the number of lines from $V(I_1)$ intersecting $C_{-1}$ twice is $\frac{(m'-1)(n'-1)}{2}$.

    In this case, if $L:\{x=a,y=b\}$ is a straight line in $V(I_2)$, then $a=\lambda+\lambda^{-1}$ and $b=\mu+\mu^{-1}$ with $\lambda^m=\mu^n=-1$. If $L$ intersects $C_{-1}$ twice, then it follows from Proposition \ref{corte} that $\lambda^{m'}=-\mu^{n'}=\pm 1$. Since $d$ is even this is a contradiction. Hence there are no lines in $V(I_2)$ intersecting $C_{-1}$ twice and the proof is complete.
    \end{proof}

    \begin{prop}\label{56}
    If $i\neq 0,\frac{d}{2}$ the number of straight lines that intersect $C_{\root^i}$ twice is $(m'-1)n'+m'(n'-1)$.
    Moreover, the previous number can be decomposed as follows.
    $$
      \left\{
      \begin{array}{ccccl}
        \frac{(m'-1)n' + m'(n'-1)}{2} & + &
        \frac{(m'-1)n' + m'(n'-1)}{2} &
        \ \text{if}\ & \text{$d, m', n'$ odd,}\\[0.15cm]
        \frac{(m'-2)n' + m'(n'-1)}{2} & + &
        \frac{m'n' + m'(n'-1)}{2} &
        \ \text{if}\ & \text{$d, n'$ odd, $m'$ even,}\\[0.15cm]
        (m'-1)n' + m'(n'-1) & + & 0 &
        \ \text{if}\ & \text{$d$ even,}
      \end{array}\right.
    $$
    where the first (resp. second) summand corresponds to the number of straight lines in $V(I_1)$
    (resp. $V(I_2)$).
    \end{prop}
    \begin{proof}
    We will assume that $d$ is odd, $m'$ is even and $n'$ is odd; the other three cases are analogous.

    Let $L:\{x=a,y=b\}$ be a straight line from $V(I_1)$, then $a=\lambda+\lambda^{-1}$ and $b=\mu+\mu^{-1}$ with $\lambda^m=\mu^n=1$ and $a,b\neq\pm 2$. From Proposition \ref{corte} it follows ($d$ being odd) that $\lambda^{m'}=1$ or $\mu^{n'}=1$ and that $\lambda^{m'}\mu^{n'}=\root^{\pm i}$. We must now consider four cases:

    \begin{enumerate}
    \item $\lambda^{m'}=1,\ \mu^{n'}=\root^i$.\hspace{0.05cm} \qquad {\bf 3)} $\ \lambda^{m'}=\root^i,\ \mu^{n'}=1$.
    \item $\lambda^{m'}=1,\ \mu^{n'}=\root^{-i}$.\qquad {\bf 4)} $\ \lambda^{m'}=\root^{-i},\ \mu^{n'}=1$.
    \end{enumerate}

    Since we are only interested in finding solutions for $(a,b)$ it is clear that we must only take into account cases 1) and 3). Now, case 1) provides $\frac{m'-2}{2}n'$ solutions in $(a,b)$ while case 3) provides $m'\frac{n'-1}{2}$. Thus, the number of straight lines in $V(I_1)$ that intersect $C_{\root^i}$ twice is $\frac{(m'-2)n'+m'(n'-1)}{2}$.

    On the other hand, let $L:\{x=a,y=b\}$ be a straight line in $V(I_2)$. Then $a=\lambda+\lambda^{-1}$ and $b=\mu+\mu^{-1}$ with $\lambda^m=\mu^n=-1$ and $a,b\neq\pm 2$. From Proposition \ref{corte} it follows ($d$ being odd) that $\lambda^{m'}=-1$ or $\mu^{n'}=-1$ and that $\lambda^{m'}\mu^{n'}=\root^{\pm i}$. Again we must consider four cases:

    \begin{enumerate}
    \item $\lambda^{m'}=-1,\ \mu^{n'}=\root^i$.\hspace{0.05cm} \qquad {\bf 3)} $\ \lambda^{m'}=\root^i,\ \mu^{n'}=-1$.
    \item $\lambda^{m'}=-1,\ \mu^{n'}=\root^{-i}$.\qquad {\bf 4)} $\ \lambda^{m'}=\root^{-i},\ \mu^{n'}=-1$.
    \end{enumerate}

    And we must only take cases 1) and 3) into account. Case 1) provides $\frac{m'}{2}n'$ solutions for $(a,b)$ while case 3) provides $m'\frac{n'-1}{2}$. Thus, the number of straight lines in $V(I_2)$ that intersect $C_{\root^i}$ twice is $\frac{m'n'+m'(n'-1)}{2}$ and the proposition follows.
    \end{proof}

  \subsection{Straight Lines Intersecting Different Components}

Now we will compute the number of lines of $\mathcal{L}$ intersecting $\mathcal{C}$ in exactly two
components. These pairs of components are $\{C_{\root^i},C_{\root^j}\}$ with $i\neq j$. We will say that a line $L$ intersects the pair $\{C_{\root^i},C_{\root^j}\}$ if $L$ intersects $C_{\root^i}$ in one point and $C_{\root^j}$ in another one.

The following elementary result will be useful later.

\begin{lemma}\label{nuLemma}
Let $\nu$ be the number of complex solutions of the system $\{\omega^2 = \root^i, \omega^d=\pm 1\}$ under
the equivalence relation $\omega_1\sim\omega_2 \Longleftrightarrow \omega_1=\omega_2^{-1}$. Then
$$
  \nu =
  \begin{cases}
  1 & \text{if}\quad i = \displaystyle\frac{d}{2},\\
  2 & \text{otherwise}.
  \end{cases}
$$
\end{lemma}

We would like to compute the number of straight lines of $\mathcal{L}$ intersecting $\{C_{\root^i},C_{\root^j}\}$
with $i\neq j$. Let $L:\{x=a, y=b\}$ be a straight line of $\mathcal{L}$ and take
$\lambda,\mu$ such that $a=\lambda+\lambda^{-1}$ and $b=\mu+\mu^{-1}$.
Then from Proposition \ref{corte} our problem is equivalent to finding the number of solutions in $(a,b)$ of the ``system''
$\{C_{\lambda^{m'}\mu^{n'}}, C_{\lambda^{m'}\mu^{-n'}}\} = \{C_{\root^i},C_{\root^j}\}$ with $i\neq j$. We
have to consider eight cases:
\begin{enumerate}
  \item $\lambda^{m'} \mu^{n'} =\root^{i},\ \lambda^{m'}\mu^{-n'}=\root^{j}$.
    \hspace{0.3cm} \qquad {\bf 5)} $\ \lambda^{m'} \mu^{n'} =\root^{i},\ \lambda^{m'}\mu^{-n'}=\root^{-j}$.
  \item $\lambda^{m'} \mu^{n'} =\root^{j},\ \lambda^{m'}\mu^{-n'}=\root^{i}$.
    \hspace{0.3cm} \qquad {\bf 6)} $\ \lambda^{m'} \mu^{n'} =\root^{-j},\ \lambda^{m'}\mu^{-n'}=\root^{i}$.
  \item $\lambda^{m'} \mu^{n'} =\root^{-i},\ \lambda^{m'}\mu^{-n'}=\root^{-j}$.
    \qquad {\bf 7)} $\ \lambda^{m'} \mu^{n'} =\root^{-i},\ \lambda^{m'}\mu^{-n'}=\root^{j}$.
  \item $\lambda^{m'} \mu^{n'} =\root^{-j},\ \lambda^{m'}\mu^{-n'}=\root^{-i}$.
    \qquad {\bf 8)} $\ \lambda^{m'} \mu^{n'} =\root^{j},\ \lambda^{m'}\mu^{-n'}=\root^{-i}$.
\end{enumerate}

Like in the proof of Proposition \ref{56} we only have to take cases 1) and 5) into account, since the other ones
do not provide more solutions for $(a,b)$. Moreover, if $i$ or $j$ equals $0$ or $\frac{d}{2}$ then
only the case 1) has to be considered.

\begin{prop}
The number of straight lines of $\mathcal{L}$ which intersect $\{C_{\root^i},C_{\root^j}\}$ with
$i\neq j$ at two different points is
$$
\begin{cases}
m'n' & \text{if}\quad (i,j)=(0,\frac{d}{2}),\\
2m'n' & \text{if}\quad i\neq 0,\frac{d}{2},\ j=0,\\
2m'n' & \text{if}\quad i\neq 0,\frac{d}{2},\ j=\frac{d}{2},\\
4m'n' & \text{if}\quad i,j\neq 0,\frac{d}{2}.
\end{cases}
$$
\end{prop}

\begin{proof}
Let us first assume that $i\neq 0,\frac{d}{2}$ and $j=0$. Then from the above
considerations we have that $\lambda^{m'}\mu^{n'} = \root^i$, $\lambda^{m'}\mu^{-n'}=1$ and
thus $\lambda^{2m'}=\root^i$, $\mu^{n'}=\lambda^{m'}$. From Lemma \ref{nuLemma} the number of
complex solutions for $(a,b)$ of the previous system is clearly $\nu m'n' = 2m'n'$.

The first and third cases are completely analogous. For the last case $i,j\neq 0,\frac{d}{2}$,
$i\neq j$ we have take into account two different systems
$\{\lambda^{m'}\mu^{n'} = \root^i$, $\lambda^{m'}\mu^{-n'}=\root^j\}$,
$\{\lambda^{m'}\mu^{n'} = \root^i$, $\lambda^{m'}\mu^{-n'}=\root^{-j}\}$.
Using again Lemma \ref{nuLemma}, one can see that both systems have
$\nu m'n'=2m'n'$ solutions and the proof is complete.
\end{proof}

  \subsection{Combinatorial Structure of $X(G)$}

  Now we are ready to present the main result of this paper that allow one to describe the
  combinatorial structure of the character variety as a function of $m$ and $n$. Recall
  that $d$ is the greatest common divisor of $m, n$ and $m = m'd$, $n=n'd$. Note
  that the proof is obvious after the above discussion.

  \begin{theo}
  The number of straight lines of $\mathcal{L}$ can be decomposed as follows.\\[0.2cm]
    {\bf a)} If $d$ is odd,
    \begin{small}
    $$
      \frac{(m'-1)(n'-1)}{2} + \frac{d-1}{2}\cdot \Big[(m'-1)n' + m'(n'-1)\Big] +
      \frac{d-1}{2}\cdot 2m'n' + \binom{(d-1)/2}{2}\cdot 4m'n',
    $$
    \end{small}
    \!\!\!\! where the first (resp. second) summand corresponds to the number of straight lines intersecting
    $C_1$ (resp. $C_{\root^i}$ with $i\neq 0$) twice and the third (resp. fourth) one corresponds to the number of
    straight lines intersecting $\{C_1, C_{\root^i}\}$ with $i\neq 0$ (resp. $\{C_{\root^i}, C_{\root^j}\}$ with
    $i,j\neq 0$, $i\neq j$) at two different points.\\[0.2cm]
    {\bf b)} If $d$ is even,
    \begin{small}
    \begin{eqnarray*}
      \frac{(m'-1)(n'-1)}{2} + \frac{(m'-1)(n'-1)}{2} +
      \frac{d-2}{2}\cdot \Big[ (m'-1)n' + m'(n'-1) \Big] +\\
      +\ m'n' + \frac{d-2}{2}\cdot 2m'n' + \frac{d-2}{2}\cdot 2m'n' + \binom{(d-2)/2}{2}\cdot 4m'n',
    \end{eqnarray*}
    \end{small}
    \!\!\!\! where the first (resp. second; third) summand corresponds to the number of straight lines
    intersecting $C_1$ (resp. $C_{-1}$; $C_{\root^i}$ with $i\neq 0,\frac{d}{2}$) twice
    and the fourth (resp. fifth; sixth; seventh) one corresponds to the number of straight lines intersecting
    $\{C_1,C_{-1}\}$ (resp. $\{C_1,C_{\root^i}\}$ with $i\neq 0,\frac{d}{2}$;\ $\{C_{-1},C_{\root^i}\}$
    with $i\neq 0,\frac{d}{2}$;\ $\{C_{\root^i},C_{\root^j}\}$ with $i,j\neq 0,\frac{d}{2}$, $i\neq j$)
    at two different points.
  \end{theo}

  \begin{ex}
  Let us use the previous theorem to describe the combinatorial structure of $X(G_{42,30})$.
  First, note that if $a|m$ and $b|n$ then
  there exists a surjective group homomorphism $G_{m,n}\twoheadrightarrow G_{a,b}$
  induced by the identity on $F_2$. Using Corollary \ref{inducedXG} one can construct
  an injective polynomial map $X(G_{a,b}) \lhook\joinrel\longrightarrow X(G_{m,n})$ coming from
  the identity on $\C^3$. Thus we have the following embeddings $X(G_{7,5}) \subset X(G_{21,15}) \subset X(G_{42,30})$
  and $X(G_{7,5}) \subset X(G_{14,10}) \subset X(G_{42,30})$. In Figure \ref{figureI} we show
  the combinatorial structure of $X(G_{42,30})$ and how $X(G_{7,5})$, $X(G_{14,10})$ and $X(G_{21,15})$
  are embedded in it.

  \begin{figure}[h]
  \centering
  \caption{Combinatorial structure of $X(G_{7,5})$, $X(G_{21,15})$, $X(G_{14,10})$ and $X(G_{42,30})$.}
  \label{figureI}
  \end{figure}
  \end{ex}

  \section{The Variety $X(G)$ as an Invariant of $G$}
\label{inv}
Once we have given a complete description of the combinatorial structure of $X(G)$, it naturally arises the question of how much information about $m$ and $n$ is codified in this structure. Clearly, if $d=1$, it will not be possible to recover any information from the combinatorial structure of $X(G)$, since it only consists of one irreducible component $V(I_3)$ cut by ${\frac{(m-1)(n-1)}{2}}$ disjoint straight lines and this latter value does not determine $m$ and $n$ even if they are coprime. For instance $(7,2)$ and $(4,3)$ have the same combinatorial structure.
Nevertheless, for $d\neq 1$ and with only one exception, we will show that it is possible to recover the value of $m$ and $n$ only from the combinatorial structure of $X(G)$.

\begin{figure}[h]
\centering
\caption{Combinatorial structure of $(7,2)$ and $(4,3)$.}
\end{figure}

Firstly, let us suppose that every irreducible component of $X(G)$ contains exactly 2 singularities. Then, by the previous section, it is easy to see that the possible values for $(m,n)$ are $(3,2)$, $(3,3)$, and $(4,2)$. Since the number of irreducible components of $X(G_{3,2})$ is 2, while $X(G_{3,3})$ and $X(G_{4,2})$ both of them contain 4 irreducible components, we can recognize the case $(m,n)=(3,2)$. Unfortunately the combinatorial structure of $X(G_{3,3})$ and $X(G_{4,2})$ are identical and it does not allow us to distinguish these cases.

\begin{figure}[h]
\centering
\caption{Combinatorial structure of $(3,3)$ and $(4,2)$.}
\end{figure}

Now, we must study the case when some of the irreducible components of $X(G)$ contain a number of singular points different from 2. These components must be precisely those coming from $V(I_3)$; let us denote them by $A_1,\dots,A_k$. Of course, if $k=1$, then $d=1$ and we have no possibility of obtaining $m$ and $n$. On the other hand, if $k\geq2$, let us denote by $a_{ij}$ the number of straight lines that intersect only $A_i$ and $A_j$ and put $M=\big( a_{ij}\big)_{1\leq i,j\leq k}$ which is a $k\times k$ symmetric matrix over the integers. Note that, for every $i\in\{1,\dots,k\}$, the sum $s_i=2a_{ii}+
{\sum_{i\neq j=1}^k a_{ij}}$ is the number of singular points of $A_i$ and consequently $s_i\neq2$ for all $i$.
By the previous section, and after a permutation of rows and columns (which preserves the set of diagonal elements of $M$), we can suppose that

$$M={\tiny\begin{pmatrix}
\frac{(m'-1)(n'-1)}{2} & 2m'n' & 2m'n' & \dots & 2m'n'\\
 & (m'-1)n'+m'(n'-1) & 4m'n' & \dots & 4m'n'\\
 & & (m'-1)n'+m'(n'-1) & \dots & 4m'n'\\
 & & &\ddots  & \vdots\\
  & & & & (m'-1)n'+m'(n'-1)
\end{pmatrix}}\ \textrm{if $d$ is odd},$$

$$M={\tiny\begin{pmatrix}
\frac{(m'-1)(n'-1)}{2} & m'n' & 2m'n' & \dots & 2m'n'\\
 & \frac{(m'-1)(n'-1)}{2} & 2m'n' & \dots & 2m'n'\\
 & & (m'-1)n'+m'(n'-1) & \dots & 4m'n'\\
 & & &\ddots  & \vdots\\
  & & & & (m'-1)n'+m'(n'-1)
\end{pmatrix}}\ \textrm{if $d$ is even}.$$

With these matrices in mind and recalling the previous section we can perform the following analysis:
\begin{itemize}
\item If $\tr M=0$, then at least $n'=1$ and three cases can occur:
\begin{enumerate}
\item If $k=2$, then $a_{12}\neq 2$ and it follows that $d=2$, $m'=a_{12}$ and $(m,n)=(2a_{12},2)$.
\item If $k>2$ and $a_{ij} = 1$ for some $i\neq j$, then  $m'=n'=1$ and $d = 2k-2$
\item If $k>2$ and $a_{ij} \neq 1$ for every $i \neq j$ then again $m'=n'=1$, but $d = 2k-1$ and we are done.
\end{enumerate}
\item If $\tr M\neq 0$, then define $a=\min\{a_{ii}\ |\ 1\leq i\leq k\}$. Two cases can occur:
\begin{enumerate}
\item If $a$ appears only once in the diagonal of $M$, then $d=2k-1$, $a=\displaystyle{\frac{(m'-1)(n'-1)}{2}}$ and $\min\{ a_{ij}\ |\ i\neq j\}=2m'n'$ so we obtain $m$ and $n$ by elementary methods.
\item If $a$ appears twice in the diagonal of $M$, then $d=2k-2$, $a=\displaystyle{\frac{(m'-1)(n'-1)}{2}}$ and $\min\{ a_{ij}\ |\ i\neq j\}=m'n'$ and we again recover the values of $m$ and $n$.
\end{enumerate}
\end{itemize}

\begin{remark}
We have seen that given the combinatorial structure of $X(G_{m,n})$ we can recover the values of $m$ and $n$ in most cases. Unfortunately there are situations where the combinatorial structure is not enough to obtain the values of $m$ and $n$; in other words, the combinatorial structure of $X(G)$ is not a complete invariant. Nevertheless in such situation we will always have only a finite number of possibilities for the pair $(m,n)$. This is clear in the case $(3,3)$ and $(4,2)$, but it is also true if $d=1$. We will give an example.

 Let us suppose that we are given the combinatorial structure of $X(G)$ and we find that $d=1$ while the number of straight lines obtained is 18. Then we know that $(m-1)(n-1)=36$ and consequently the only possibilities ($m$ and $n$ being coprime) are $(37,2)$, $(19,3)$ and $(13,4)$. This situation is general by virtue of the prime decomposition. However, even if $d=1$ there are cases where we can recover $m$ and $n$. For instance if there are $p$ straight lines where $p$ is a prime of the form $3k+2$ (there are infinetely many) the only possible solution for $(m,n)$ is $(2p+1,2)$.

 It would be interesting to find or construct another invariant from the character variety which could be shown to be complete.
\end{remark}

    \section{Irreducible Components of $R(G_{m,n})$}
    \label{irrR}
    The main goal of this section will be to compute the number of irreducible components of $R(G_{m,n})$ and to study the behavior of the projection $t:R(G_{m,n})\longrightarrow X(G_{m,n})$. In particular we will see that like in \cite[Rem. 3.17.]{POR} (where the variety of characters in $PSL(2,\C)$ is considered), the projection $t$ always induces a bijection between irreducible components.

    Given the set $R(G_{m,n})$ we define the following subsets:
    $$Irr=\{\rho\in R(G_{m,n})\ |\ \rho\ \textrm{is irreducible}\}$$
    $$Red=\{\rho\in R(G_{m,n})\ |\ \rho\ \textrm{is reducible}\}$$
    $$Ab=\{\rho\in R(G_{m,n})\ |\ \rho(G_{m,n})\ \textrm{is abelian}\}$$
    $$M_{R}=\{\rho\in Red\ |\ \rho(G_{m,n})\ \textrm{is metabelian}\}$$
    Clearly we have the partitions $R(G_{m,n})=Irr\cup Red=Irr\cup Ab\cup M_R$. We know that $Ab$ is an algebraic set, while $Irr$ and $M_R$ are not. We can obtain a decomposition of $R(G_{m,n})$ into closed subsets just by taking closures, but the unions will no longer be disjoint. Namely, we have $R(G_{m,n})=\overline{Irr}\cup Ab\cup\overline{M_R}$. Moreover, since it can be seen that $\overline{M_R}\subseteq\overline{Irr}$ we have that $R(G_{m,n})=\overline{Irr}\cup Ab$ and we will study these two subsets separatedly in order to count the number of irreducible components of $R(G_{m,n})$.

\vspace{0.25cm}

    Note that we have $t(\overline{Irr})=V(I_1)\cup V(I_2)$, $t(Ab)=V(I_3)=t(Red)$ and $t(\overline{M_R})=\left(V(I_1)\cup V(I_2)\right)\cap V(I_3)=t(\overline{Irr})\cap t(Ab)$. In particular this implies ($t$ being continuous) that the previous decomposition is not redundant.

    \subsection{The Irreducible Components of $\overline{Irr}$}
    The following result can be shown by elementary facts of general topology.
    \begin{lemma}
    Let $\varphi : X \longrightarrow Y$ be a continuous surjective map between two
    topological spaces. Then the number of irreducible components of $Y$ is
    less than or equal to the number of irreducible components of $X$. Moreover,
    if this numbers are different then there exist $X_1$ and $X_2$ two different components of $X$
    such that $\overline{\varphi(X_1)} \subseteq \overline{\varphi(X_2)} $.
    \end{lemma}

    Note that if $X_1$ and $X_2$ are two different components of $R(G_{m,n})$ which contain
    irreducible representations, then $t(X_1)\cap t(X_2)=\emptyset$ and thus
    $t|_{\overline{Irr}}: \overline{Irr}\longrightarrow
    t(\overline{Irr})$ preserves the number of irreducible components.

    \begin{remark}
    Let $G$ be an arbitrary finitely presented group. From Propositions 1.5.2 and 1.1.1 in \cite{CS83},
    $t(X_1)\cap t(X_2) \subseteq t(Red)$ and hence $t(X_1)\nsubseteq t(X_2)$.
    Therefore $t|_{\overline{Irr}}: \overline{Irr}\longrightarrow
    t(\overline{Irr})$ always preserves the number of irreducible.
    \end{remark}

    As we proved in the previous section,
    $t(\overline{Irr})$ has either $\frac{(|m|-1)(|n|-1)}{2}$ irreducible components if $d$ is odd or
    $\frac{(|m|-1)(|n|-1)+1}{2}$ if $d$ is even. Thus we have found the number of irreducible components of
    $\overline{Irr}$. Note that that due to \cite[Cor. 1.5.3.]{CS83} all these irreducible
    components are of dimension 4.

   \medskip
    Since it is known that $\dim R(G_{m,n})=4$ and we will see that $\dim Ab=\dim \overline{M_R}=3$,
    we have; in particular, the following result:


    \begin{theo}
    The number or irreducible 4-dimensional components of $R(G_{m,n})$ is:
    \begin{itemize}
    \item[a)] $\displaystyle{\frac{(|m|-1)(|n|-1)}{2}}$ if $d$ is odd.
    \item[b)] $\displaystyle{\frac{(|m|-1)(|n|-1)+1}{2}}$ if $d$ is even.
    \end{itemize}
    \end{theo}

    This result can be found in \cite[Theorem A]{LI2}, where a more direct approach is used.
    Note that we obtained it as a consequence of our study of $X(G_{m,n})$.

    \subsection{The Irreducible Components of $Ab$}
    This section is devoted to count the number of irreducible components of $Ab$. First we note that
    $Ab=R(G_{m,n}^{\textrm{ab}})$ where $G^{\textrm{ab}}=G/G'$ denotes the abelianization of $G$. In our case
    $G_{m,n}^{\textrm{ab}}=\langle x,y\ |\ x^m=y^n,\ [x,y]=1\rangle$. In the following lemma we give another
    presentation of $G_{m,n}^{\textrm{ab}}$ which will be easier to work with.
    \begin{lemma}
    $G_{m,n}^{\textrm{ab}}\cong H_d=\langle a,b\ |\ a^d=1=[a,b]\rangle$, where $d=\gcd(m,n)$.
    \end{lemma}
    \begin{proof}
    Put $m'=\frac{m}{d}$ and $n'=\frac{n}{d}$ and consider Bezout's identity $\alpha m-\beta n=d$. The claimed
    isomorphism is then given by $\phi:G^{\textrm{ab}}\longrightarrow H_d$ with $\phi(x)=b^{n'}a^{\alpha}$,
    $\phi(y)=b^{m'}a^{\beta}$ and $\psi:H_d\longrightarrow G^{\textrm{ab}}$ with $\psi(a)=x^{m'}y^{-n'}$,
    $\psi(b)=x^{-\beta}y^{\alpha}$.
    \end{proof}

    Thus, $Ab\cong R(H_d)$ and we will count the irreducible components of the latter. To do so we
    introduce some notation. If $\xi\in\{1,\dots,\root^{d-1}\}$ is a $d$-th root of unity
    we put $A_{\xi}=\left(\begin{smallmatrix} \xi & 0\\ 0 & \xi^{-1}\end{smallmatrix}\right)$
    and we define the set $V_{\xi}=\{P^{-1}A_{\xi}P\ |\ P\in SL(2,\C)\}$. Note that $V_1=\{I_2\}$,
    $V_{-1}=\{-I_2\}$ and if $\xi\neq\pm1$ then $V_{\xi}$ is an irreducible affine algebraic variety
    of dimension 2 (see \cite[Cor. 1.5.]{LI1}).

    \begin{lemma}
    If $X\in SL(2,\C)$ is such that $X^d=I_2$, then $X\in V_{\xi}$ for some $\xi$ $d$-th
    root of unity.
    \end{lemma}
    \begin{proof}
    Since $\C$ is algebraically closed there exists $P\in SL(2,\C)$ such that
    $PXP^{-1}=\left(\begin{smallmatrix} a & \alpha\\ 0 & a^{-1}\end{smallmatrix}\right) = Y$.
    Clearly $I_2=Y^d=\left(\begin{smallmatrix} a^d & \alpha h_d(a+a^{-1})\\ 0 & a^{-d}\end{smallmatrix}\right)$
    so $a^d=1$ and $\alpha h_d(a+a^{-1})=0$ and two cases arise.
    \begin{enumerate}
    \item If $a=\pm 1$, since $h_d(\pm 2)\neq0$ it must be $\alpha=0$ and $X=Y=\pm I_2\in V_{\pm1}$.
    \item If $a=\xi\neq\pm 1$ then $a\neq a^{-1}$ and $X$ is diagonalizable so there exists
    $P\in SL(2,\C)$ such that $X=P^{-1}A_{\xi}P\in V_{\xi}$.
    \end{enumerate}
    \end{proof}

    If we now define $M_{\xi}=\{(A,B)\in SL(2,\C)\ |\ A\in V_{\xi},\ [A,B]=I_2\}$, then the previous lemma
    shows that $R(H_d)=\bigcup_{i=0}^{d-1}M_{\root^i}$. Clearly $M_1=\{I_2\}\times SL(2,\C)$ and
    $M_{-1}=\{-I_2\}\times SL(2,\C)$ are irreducible affine algebraic varieties of dimension 3.
    Now, we want to study the case $M_{\xi}$ with $\xi\neq\pm 1$.

    \begin{prop}
    $M_{\xi}$ is an affine irreducible algebraic variety of dimension 3
    for all $\xi$ $d$-th root of unity.
    \end{prop}
    \begin{proof}
    We can assume $\xi\neq\pm1$. In this case we define $\Psi:M_{\xi}\longrightarrow V_{\xi}
    \times \C^{*}$ as follows: given $(A,B)\in M_{\xi}$ there exists $P\in SL(2,\C)$ such that
    $PAP^{-1}=A_{\xi}$, now since $A$ and $B$ commute, $PBP^{-1}$ and $A_{\xi}$ must also commute
    and it follows that $PBP^{-1}=\left(\begin{smallmatrix} b & 0\\ 0 & b^{-1}\end{smallmatrix}\right)$
    must be diagonal. We define $\Psi(A,B)=(A,b)$.

    Let us see that $\Psi$ is well defined: if $PAP^{-1}=A_{\xi}=QAQ^{-1}$, then $QP^{-1}$ commute with
    $A_{\xi}$ and it must be diagonal. Consequently $QBQ^{-1}=QP^{-1}\left(\begin{smallmatrix} b & 0\\ 0 &
    b^{-1}\end{smallmatrix}\right) (QP^{-1})^{-1}=\left(\begin{smallmatrix} b & 0\\ 0 & b^{-1}
    \end{smallmatrix}\right)$ and $\Psi(A,B)$ does not depend on the choice of $P$.

    Now, we claim that $\Psi$ is bijective. Let us suppose that $\Psi(A_1,B_1)=(A_1,b_1)=(A_2,b_2)=
    \Psi(A_2,B_2)$, then $A_1=A_2$ and $b_1=b_2$. Since $A_1=A_2$ and $\Psi$ is well defined there
    must exist $P\in SL(2,\C)$ such that $PB_1P^{-1}=\left(\begin{smallmatrix} b_1 & 0\\ 0 &
    b_1^{-1}\end{smallmatrix}\right)=\left(\begin{smallmatrix} b_2 & 0\\ 0 & b_2^{-1}\end{smallmatrix}\right)
    =PB_2P^{-1}$ so we have that $(A_1,B_1)=(A_2,B_2)$ and $\Psi$ is injective. Since surjectivity of
    $\Psi$ is obvious the claim follows.

    To finish the proof it is enough to observe that $\Psi$ induces a birational equivalence between $M_{\xi}$ and
    $V_{\xi}\times\C$.
    \end{proof}

    Now, in order to be able to count the number of irreducible components, we must remove the redundant
    components in the decomposition $R(H_d)=\bigcup_{i=0}^{d-1} M_{\root^i}$. This is done as follows.

    \begin{prop}
    $M_{\xi}=M_{\eta}$ if and only if $\xi=\eta^{\pm1}$. Moreover, if
    $M_{\xi}\neq M_{\eta}$, then they are disjoint.
    \end{prop}
    \begin{proof}
    $M_{\xi}=M_{\eta}\Leftrightarrow V_{\xi}=V_{\eta}\Leftrightarrow A_{\xi}\in
    V_{\eta}\stackrel{\text{\cite[1.4.]{LI1}}}{\Leftrightarrow} \tr A_{\xi}=\tr A_{\eta}
    \Leftrightarrow \xi=\eta^{\pm1}$.
    \end{proof}

    As a consequence of this proposition we have that $Ab$ has $\frac{d+1}{2}$ irreducible components
    if $d$ is odd and $\frac{d+2}{2}$ if $d$ is even.
    
    \subsection{Counting the Irreducible Components of $R(G_{m,n})$}
    We can now combine the results obtained in the previous sections to explicitly compute the number of irreducible components of $R(G_{m,n})$. Namely we have the following.

    \begin{theo}
    The number of irreducible components of $R(G_{m,n})$ is
    $$\begin{cases}
    \displaystyle{\frac{(|m|-1)(|n|-1)}{2}+\frac{d+1}{2}} & \textrm{if $d$ is odd},\\
    \ &\ \\
    \displaystyle{\frac{(|m|-1)(|n|-1)+1}{2}+\frac{d+2}{2}} & \textrm{if $d$ is even},
    \end{cases}$$
    where the first summand corresponds to the number of irreducible components in $\overline{Irr}$ and the second one to the irreducible components of $Ab$.
    \end{theo}

     Note that, in the light of Theorem \ref{compx}, the projection $t:R(G_{m,n})\longrightarrow X(G_{m,n})$ preserves the number of irreducible components. 
     In particular, if we consider the restrictions $t|_{\overline{Irr}}:\overline{Irr}\longrightarrow V(I_1)\cup V(I_2)$ and $t|_{Ab}:Ab\longrightarrow V(I_3)$, we have that both of them induce bijections between irreducible components.

     Also note that, due to \cite[Cor. 1.5.3.]{CS83}, if $Irr_0$ is an irreducible component of $\overline{Irr}$, then $\dim Irr_0=\dim t(Irr_0)+3$. In our case $\dim Irr_0=4$ and $\dim t(Irr_0)=1$.  Nevertheless, it is known that the result is not true in general and, for instance, we have that $\dim Ab_0=3$ and $\dim t(Ab_0)=1$ for every irreducible component of $Ab$.

    \subsection{The Irreducible Components of $\overline{M_R}$}
    Although we have just computed the number of irreducible components of $R(G_{m,n})$ it can be interesting to compute the number of components of $\overline{M_R}$ and to see if $t|_{\overline{M_R}}$ preserves such number. Recall that $M_R$ is the set of reducible metabelian representations of $G_{m,n}$ and $t(\overline{M_R})=t(\overline{Irr})\cap t(Ab)$. 
    In the following lemma we give some properties of such representations which will be useful
    in the sequel.

    \begin{lemma}\label{propmeta}
    Let $\rho\in M_R$. Then $\rho(x)^m=\rho(y)^n=\pm I_2$ and $\tr\rho(x), \tr\rho(y)\neq 2$.
    \end{lemma}
    \begin{proof}
    We can assume that $\rho(x)=X=\left(\begin{smallmatrix} a & \alpha\\ 0 & a^{-1}\end{smallmatrix}\right)$ and
    $\rho(y)=Y=\left(\begin{smallmatrix} b & \beta \\ 0 & b^{-1}\end{smallmatrix}\right)$.
    Let us suppose that $\tr X=2$. Then $a=1$ and $b^n=1$ and two cases arise:
    \begin{enumerate}
    \item If $b=\pm1$ then $\rho(G)$ is abelian, which contradicts the hypothesis $p\in M_R$.
    \item If $b\neq\pm1$ then $I_2=Y^n=X^m=\left(\begin{smallmatrix} 1 & m\alpha\\ 0 & 1\end{smallmatrix}\right)$
    and $m\alpha=0$, thus $\alpha=0$ and $X=I_2$ which implies again that $\rho(G)$ is abelian.
    \end{enumerate}
    Analogously it can be proved that $\tr X\neq -2$ and $\tr Y\neq\pm2$.

    Now, suppose that $X^m=Y^n\neq\pm I_2$, in which case $a^m,b^n\neq\pm1$ and
    $h_m(a+a^{-1})\neq0\neq
    h_n(b+b^{-1})$. From $A^m=B^n$ it follows that $\alpha h_m(a+a^{-1})=\beta h_n(b+b^{-1})$, so $\alpha=0$
    if and only if $\beta=0$ and, since $\alpha=\beta=0$ implies that $\rho(G)$ is abelian we deduce that
    $\alpha\beta\neq0$. Finally $a^m+a^{-m}=b^n+b^{-n}\Leftrightarrow(a-a^{-1})h_m(a+a^{-1})=
    (b+b^{-1})h_n(b+b^{-1})\Leftrightarrow \beta(a-a^{-1})=\alpha(b-b^{-1})\Leftrightarrow$
    $\rho(G)$ is abelian. This, again, is a contradiction and the lemma follows.
    \end{proof}

    We will now introduce some notation.
    $$\Theta=\{\xi\ |\ \xi^m=\pm1,\ \xi\neq\pm1\}=\{\xi\ |\ \xi^m=1,\ \xi\neq\pm1\}\cup\{\xi\ |\
    \xi^m=-1,\ \xi\neq-1\}=\Theta^{+}\cup\Theta^{-}$$
    $$\Upsilon=\{\eta\ |\ \eta^n=\pm1,\ \eta\neq\pm1\}=\{\eta\ |\ \eta^n=1,\ \eta\neq\pm1\}\cup\{\eta\ |\
    \eta^n=-1,\ \eta\neq-1\}=\Upsilon^{+}\cup\Upsilon^{-}$$
    Now, given $\xi\in\Theta$ and $\eta\in\Upsilon$ we put $A_{\xi}=\left(\begin{smallmatrix} \xi & 0\\ 0 &
    \xi^{-1}\end{smallmatrix}\right)$ and $B_{\eta}=\left(\begin{smallmatrix} \eta & 1 \\ 0 &
    \eta^{-1}\end{smallmatrix}\right)$. Finally, let us define $V_{\xi,\eta}=\{(P^{-1}A_{\xi}P,P^{-1}
    B_{\eta}P)\ |\ P\in SL(2,\C)\}$.

    \begin{lemma}
    If $\rho\in M_R$, there exist $\xi\in\Theta$ and $\eta\in\Upsilon$ such that $(\rho(x),\rho(y))\in
    V_{\xi,\eta}$.
    \end{lemma}
    \begin{proof}
    Put $X=\rho(x)$ and $Y=\rho(y)$. By Lemma \ref{propmeta} we have that there exists $P\in SL(2,\C)$ such that
    $PXP^{-1}=\left(\begin{smallmatrix} \xi & \alpha\\ 0 & \xi^{-1}\end{smallmatrix}\right)=X'$
    and $PYP^{-1}=\left(
    \begin{smallmatrix} \eta & \beta \\ 0 & \eta^{-1}\end{smallmatrix}\right)=Y'$ for some $\xi\in\Theta$,
    $\eta\in\Upsilon$ and $\alpha,\beta\in\C$. By straightforward computations it is easy to see that there
    exists an upper triangular matrix $Q\in SL(2,\C)$ such that $Q^{-1}X'Q=A_{\xi}$ and $Q^{-1}Y'Q=B_{\eta}$.
    This completes the proof.
    \end{proof}

    As a consequence of this lemma we obtain a decomposition
    $$\overline{M_R}=\bigcup_{\substack{\xi\in\Theta^{+}\\
    \eta\in\Upsilon^{+}}}\overline{V_{\xi,\eta}}\cup\bigcup_{\substack{\xi\in\Theta^{-}\\
    \eta\in\Upsilon^{-}}}\overline{V_{\xi,\eta}}.$$
    Therefore, in order to find the dimension and the number of irreducible components of $\overline{M_R}$ we will study
    each $\overline{V_{\xi,\eta}}$ separatedly.

    \begin{prop}
    Given $\xi\in\Theta^{\pm}$ and $\eta\in\Upsilon^{\pm}$, the set $\overline{V_{\xi,\eta}}$ is an affine
    irreducible algebraic variety of dimension 3.
    \end{prop}
    \begin{proof}
    Let us define $\Phi:V_{\xi,\eta}\longrightarrow PSL(2,\C)$ as follows: given $(A,B)\in
    V_{\xi,\eta}$, there exists  $P\in SL(2,\C)$ such that $PAP^{-1}=A_{\xi}$ and $PBP^{-1}=B_{\eta}$;
    we define $\Phi(A,B)=[P]$.

    Now we will see that $\Phi$ is well defined. If $PAP^{-1}=A_{\root}=QAQ^{-1}$, and $PBP^{-1}=B_{\eta}=QBQ^{-1}$ it follows that $PQ^{-1}$ commutes with $A_{\xi}$, so it must be diagonal. Moreover as $PQ^{-1}$ commutes with $B_{\eta}$ and it is diagonal it must be $PQ^{-1}=\pm I_2$ so $P=\pm Q$ and $[P]=[Q]$ in $PSL(2,\C)$.

    Since $\Phi$ is trivially surjective, it is enough to prove the injectivity. If $(A_1,B_1),(A_2,B_2)\in V_{\xi,\eta}$ are such that $\Phi(A_1,B_1)=\Phi(A_2,B_2)=[P]$, then either $PA_1P^{-1}=A_{\xi}=PA_2P^{-1}$ or $PA_1P^{-1}=(-P)A_2(-P)^{-1}$. But in any case
    $A_1=A_2$ and analogously $B_1=B_2$.

    Now, since $\Phi$ clearly induces a birational equivalence between $V_{\xi,\eta}$ and
    $PSL(2,\C)\cong SO_3(\C)$, the proof is complete.
    \end{proof}

    Note that $V_{\xi_i,\eta_j}=V_{\xi_k,\eta_l}$ if and only if $\xi_i=\xi_k$ and $\eta_j=\eta_l$. Moreover if $V_{\xi_i,\eta_j}\neq V_{\xi_k,\eta_l}$ then they are disjoint. Finally, as in the previous sections, we have to remove the redundant components in our decomposition of $\overline{M_R}$.

    \begin{prop}
    $\overline{V_{\xi_i,\eta_j}}=\overline{V_{\xi_k,\eta_l}}$ if and only if $\xi_i=\xi_k$ and $\eta_j=\eta_l$.
    \end{prop}
    \begin{proof}
    Let us suppose that $\overline{V_{\xi_i,\eta_j}}=\overline{V_{\xi_k,\eta_l}}$, then $(\xi_i+\xi_i^{-1},\eta_j+\eta_j^{-1},\xi_i\eta_j+(\xi_i\eta_j)^{-1})=t\left(\overline{V_{\xi_i,\eta_j}}\right)=t\left(\overline{V_{\xi_k,\eta_l}}\right)=(\xi_k+\xi_k^{-1},\eta_l+\eta_l^{-1},\xi_k\eta_l+(\xi_k\eta_l)^{-1})$. This implies that either $\xi_i=\xi_k$ and $\eta_j=\eta_l$ or $\xi_i=\xi_k^{-1}$ and $\eta_j=\eta_l^{-1}$.

    Now, if $\xi_i=\xi_k^{-1}$ and $\eta_j=\eta_l^{-1}$ we know that $V_{\xi_i,\eta_j}\cap V_{\xi_k,\eta_l}=\emptyset$. Consequently $V_{\xi_i,\eta_j}\subseteq\overline{V_{\xi_k,\eta_l}}-V_{\xi_k,\eta_l}$. This is a contradiction since $\dim(\overline{V_{\xi_k,\eta_l}}\setminus V_{\xi_k,\eta_l})<\dim V_{\xi_k,\eta_l}=\dim V_{\xi_i,\eta_j}$ (see \cite[\S 8.3.]{HUM} for instance).

   The converse is obvious.
    \end{proof}

    \begin{cor}
    The number of irreducible components of $\overline{M_R}$ is 
    $$\begin{cases}
    2(|m|-1)(|n|-1) & \textrm{if $d$ is odd},\\
    \ &\ \\
    2[(|m|-1)(|n|-1)+1] & \textrm{if $d$ is even}.
    \end{cases}$$
    \end{cor}
    
     It can be seen that $t|_{\overline{M_R}}$ maps every irreducible component of $\overline{M_R}$ to a single point which lies in $t(\overline{Irr})\cap t(Ab)$. Moreover, there are 2 irreducible components of $\overline{M_R}$ mapping to the same point in $t(\overline{Irr})\cap t(Ab)$, namely: $$\{V_{\xi,\eta},V_{\xi^{-1},\eta^{-1}}\}\stackrel{t}{\longrightarrow}\left(\xi+\xi^{-1},\eta+\eta^{-1},\xi\eta+(\xi\eta)^{-1}\right)$$ for all $(\xi,\eta)\in\Theta\times\Upsilon$.
     
     Also observe that although \cite[Cor. 1.5.3.]{CS83} cannot be applied in this situation and thus, if $M_{R_0}$ is an irreducible component of $\overline{M_R}$, then we have that $3=\dim M_{R_0}=\dim t(M_{R_0})+3=0+3$.

  \section{Comments and Applications}
  \label{com}

    \subsection{The Abelian Component}
    Since $V(I_3)=t(Ab)$ we will call it the abelian component of $X(G_{m,n})$. In Section \ref{descri}
    we gave a simple description of this component and its decomposition into irreducible components.
    Namely, we proved that $V(I_3)=\bigcup_{i=1}^{d}C_{\root_i}$ with $C_{\root_i}=\{(u+u^{-1},v+v^{-1},
    uv+(uv)^{-1})\ |\ v,v\in\C^{*},\ u^{m'}=\root_i v^{n'}\}$ and $C_{\root_i}=C_{\root_j}$ if and only if
    $\root_i=\root_j^{\pm1}$.

    Now, if $\mathcal{C}_{\root_i}$ is the plane curve given by $u^{m'}=\root_i v^{n'}$, it is clear that
    $\mathcal{C}_{\root_i}$ is isomorphic to $\mathcal{C}_1$ for all $i=1,\dots,d$. Since $C_{\root_i}$
    is birationally equivalent to $\mathcal{C}_{\root_i}$ it follows that $C_{\root_i}$ is birationally
    equivalent to $C_1$ for all $i=1,\dots,d$. Therefore, if we want to study the curves $C_{\root_i}$ it
    is enough to study $C_1$ and, as we will see, it admits a parametrization.

    \begin{lemma}
    Let $r,s$ be coprime integers. If $A,B$ are two commuting matrices in $SL(2,\C)$ such that $A^r=B^s$,
    then there exists $C\in SL(2,\C)$ such that $A=C^s$ and $B=C^r$.
    \end{lemma}

    \begin{proof}
    It is a straightforward consequence of the fact that the group $\langle a,b\ |\ a^r=b^s,\ [a,b]=1\rangle$
    is cyclic if $\gcd(r,s)=1$.
    \end{proof}

    With this lemma we are ready to give a parametrization of $C_1$.

    \begin{prop}
    $C_1=\{(f_{n'}(t),f_{m'}(t),f_{n'+m'}(t))\ |\ t\in\C\}$.
    \end{prop}
    \begin{proof}
    We know that $C_1=\{(u+u^{-1},v+v^{-1},uv+(uv)^{-1})\ |\ v,v\in\C^{*},\ u^{m'}=v^{n'}\}$ with
    $\gcd(m',n')=1$. Now if we put $A=\left(\begin{smallmatrix} u & 0\\ 0 & u^{-1}\end{smallmatrix}\right)$
    and $B=\left(\begin{smallmatrix} v & 0\\ 0 & v^{-1}\end{smallmatrix}\right)$, we can apply the previous
    lemma to find $C\in SL(2,\C)$ such that $A=C^{n'}$ and $B=C^{m'}$. Thus, $u+u^{-1}=\tr A=\tr C^{n'} =
    f_{n'}(\tr C)$; analogously $v+v^{-1}=f_{m'}(\tr C)$ and $uv+(uv)^{-1}=f_{m'+n'}(\tr C)$. Finally it is
    enough to put $t=\tr C$ and the result follows.
    \end{proof}

    Finally we want to point out that if we project $C_1$ over the plane $Z=0$ we essentially obtain a
    Lissajous figure. This is due to the relation between $\{f_k\}_{k\in\mathbb{N}}$ and the Chebyshev
    polynomials of the first kind $\{T_k\}_{k\in\mathbb{N}}$ and the well-known parametrization of the
    Lissajous curves.

    \subsection{The Mirror Image of $K_{m,n}$ and its Character Variety}
    In this section we are interested in the relationship between the character variety of the
    torus knot $K_{m,n}$ and that of its mirror image (which is precisely $K_{m,-n}$). We know that
    $X(G_{m,n})=V(J_{m,n})=V(I_1^{m,n})\cup V(I_2^{m,n})\cup V(I_3^{m,n})$ and
    $X(G_{m,-n})=V(J_{m,-n})=V(I_1^{m,-n})\cup V(I_2^{m,-n})\cup V(I_3^{m,-n})$.

    Let us define $\phi:R(G_{m,n})\longrightarrow R(G_{m,-n})$ given by $\phi(A,B)=(A,B^{-1})$.
    Clearly it is an isomorphism and due to Lemma \ref{rax} it induces an isomorphism (an involution in fact)
    $\psi:X(G_{m,n})\longrightarrow X(G_{m,-n})$ which can be seen to be given by $\psi(X,Y,Z)=(X,Y,XY-Z)$.
    By definition we have that $\psi(V(I_i^{m,n}))=V(I_i^{m,-n})$ for $i=1,2,3$. Now it is easy to see that,
    due to Remark \ref{signok}., $\psi$ fixes (not pointwise) $V(I_1^{m,n})$ and $V(I_2^{m,n})$ so $X(G_{m,n})$ and
    $X(G_{m,-n})$ share the set of straight lines. Moreover we have the decompositions
    $$V(I_3^{m,n})=\bigcup_{i=1}^d C_{\root_i}^{m,n},\quad V(I_3^{m,-n})=\bigcup_{i=1}^d C_{\root_i}^{m,-n}$$
    and it can be trivially seen that $\psi(C_{\root_i}^{m,n})=C_{\root_i}^{m,-n}$.

    We will end this section studying the intersection $\mathcal{I}=V(I_3^{m,n})\cap V(I_3^{m,-n})$.
    Recall that $V(I_3^{m,n})=\{(u+u^{-1},v+v^{-1},uv+(uv)^{-1})\ |\ u^m=v^n\}$ and
    $V(I_3^{m,-n})=\{(u+u^{-1},v+v^{-1},uv+(uv)^{-1})\ |\ u^m=v^{-n}\}$. An elementary analysis of the situation
    allows us to see that $(u+u^{-1},v+v^{-1},uv+(uv)^{-1})\in\mathcal{I}$ if and only if $u^m=v^n=\pm1$. This
    implies that $\textrm{card}(\mathcal{I})=mn+1$.

    \subsection{Another Description for the Character Variety of $K_{m,2}$}
    If $m$ is an odd integer, then $K_{m,2}$ is a two-bridge knot and it is possible to find a presentation
    of $G_{m,2}:=G(K_{m,2})$ generated by two meridians. Namely, we have an isomorphism
    $$
      H_m=\langle x,y\ |\ \overbrace{xyxy\dots yx}^{\textrm{length $m$}}=
      \underbrace{yxyx\dots xy}_{\textrm{length $m$}}\rangle\cong\langle a,b\ |\ a^m=b^2\rangle=G_{m,2}
    $$
    given by $x\mapsto b^{-1}a^{\frac{m+1}{2}}$, $y\mapsto a^{-\frac{m-1}{2}}b$ which, due to Corollary 1.6.,
    induces an isomorphism between $X(H_m)$ and $X(G_{m,2})$ given by $(X,Y,Z)\mapsto
    \left(F_{\frac{m+1}{2},1}(X,Y,Z),X\right)$
    for all $(X,Y,Z)\in X(G_{m,2})$. Now, it can be proved using the techniques in \cite{OLL} that
    $$
      X(H_m)\cong\{(X,Y)\in\C^2\ |\ (X^2-Y-2)\sigma_m(Y)=0\}
    $$
    and we have a description of $X(G_{m,2})$ as an algebraic complex curve in $\C^2$.

    Since for $|m|,|n|>2$ the knot 
    $K_{m,n}$ has more than two bridges, the Wirtinger presentation of $G_{m,n}$ has more than two generators and the
    character variety is not a curve in $\C^2$. Nevertheless it would be interesting to find an explicit isomorphism
    between the Wirtinger presentation and the one given in Section \ref{XG} in order to find an easier description of
    $X(G_{m,n})$ and, in particular, of its abelian component.

  \section*{Acknowledgments}

    We would like to thank Enrique Artal and Mar\'{\i}a Teresa Lozano for
    their constant support and motivation in our work over the years.
    Also we wish to express our gratitude to Jos\'{e} Ignacio Cogolludo for his proofreading
    and constant advise.



    The first author is partially supported by the Spanish projects MTM2007-67908-C02-01 and FQM-333,
    and the second author by the project MTM2007-67884-C04-02.


  \bibliography{./refTorusLinks}

\begin{thebibliography}{10}

\bibitem{CS83}
Marc Culler and Peter~B. Shalen.
\newblock Varieties of group representations and splittings of {$3$}-manifolds.
\newblock {\em Ann. of Math. (2)}, 117(1):109--146, 1983.

\bibitem{GM93}
F.~Gonz{\'a}lez-Acu{\~n}a and Jos{\'e}~Mar{\'{\i}}a Montesinos-Amilibia.
\newblock On the character variety of group representations in {${\rm
  SL}(2,{\bf C})$} and {${\rm PSL}(2,{\bf C})$}.
\newblock {\em Math. Z.}, 214(4):627--652, 1993.

\bibitem{Hartshorne77}
Robin Hartshorne.
\newblock {\em Algebraic geometry}.
\newblock Springer-Verlag, New York, 1977.
\newblock Graduate Texts in Mathematics, No. 52.

\bibitem{POR}
Michael Heusener and Joan Porti.
\newblock The variety of characters in {${\rm PSL}\sb 2(\Bbb C)$}.
\newblock {\em Bol. Soc. Mat. Mexicana (3)}, 10(Special Issue):221--237, 2004.

\bibitem{HLM2}
Hugh~M. Hilden, Mar{\'{\i}}a~Teresa Lozano, and Jos{\'e}
  Mar{\'{\i}}a~Montesinos Amilibia.
\newblock Character varieties and peripheral polynomials of a class of knots.
\newblock {\em J. Knot Theory Ramifications}, 12(8):1093--1130, 2003.

\bibitem{HLM92}
Hugh~M. Hilden, Mar{\'{\i}}a~Teresa Lozano, and Jos{\'e}~Mar{\'{\i}}a
  Montesinos-Amilibia.
\newblock On the character variety of group representations of a {$2$}-bridge
  link {$p/3$} into {${\rm PSL}(2,{\bf C})$}.
\newblock {\em Bol. Soc. Mat. Mexicana (2)}, 37(1-2):241--262, 1992.
\newblock Papers in honor of Jos{\'e} Adem (Spanish).

\bibitem{HLM95}
Hugh~M. Hilden, Mar{\'{\i}}a~Teresa Lozano, and Jos{\'e}~Mar{\'{\i}}a
  Montesinos-Amilibia.
\newblock On the arithmetic {$2$}-bridge knots and link orbifolds and a new
  knot invariant.
\newblock {\em J. Knot Theory Ramifications}, 4(1):81--114, 1995.

\bibitem{HLM1}
Hugh~M. Hilden, Mar{\'{\i}}a~Teresa Lozano, and Jos{\'e}~Mar{\'{\i}}a
  Montesinos-Amilibia.
\newblock On the character variety of tunnel number 1 knots.
\newblock {\em J. London Math. Soc. (2)}, 62(3):938--950, 2000.

\bibitem{HL05}
Hugh~M. Hilden, Mar{\'{\i}}a~Teresa Lozano, and Jos{\'e}~Mar{\'{\i}}a
  Montesinos-Amilibia.
\newblock Peripheral polynomials of hyperbolic knots.
\newblock {\em Topology Appl.}, 150(1-3):267--288, 2005.

\bibitem{HUM}
James~E. Humphreys.
\newblock {\em Linear algebraic groups}.
\newblock Springer-Verlag, New York, 1975.
\newblock Graduate Texts in Mathematics, No. 21.

\bibitem{LI1}
Sal Liriano.
\newblock Algebraic geometric invariants for a class of one-relator groups.
\newblock {\em J. Pure Appl. Algebra}, 132(1):105--118, 1998.

\bibitem{LI2}
Sal Liriano.
\newblock Irreducible components in an algebraic variety of representations of
  a class of one relator groups.
\newblock {\em Internat. J. Algebra Comput.}, 9(1):129--133, 1999.

\bibitem{LM85}
Alexander Lubotzky and Andy~R. Magid.
\newblock Varieties of representations of finitely generated groups.
\newblock {\em Mem. Amer. Math. Soc.}, 58(336):xi+117, 1985.

\bibitem{OLL}
Antonio~M. Oller-Marc{\'e}n.
\newblock The {${\rm SL}(2,{\bf C})$} character variety of a class of torus
  knots.
\newblock {\em Extracta Math.}, to appear.

\bibitem{RI90}
Theodore~J. Rivlin.
\newblock {\em Chebyshev polynomials}.
\newblock Pure and Applied Mathematics (New York). John Wiley \& Sons Inc., New
  York, second edition, 1990.
\newblock From approximation theory to algebra and number theory.

\bibitem{RUD}
Ze{\'e}v Rudnick.
\newblock Representation varieties of solvable groups.
\newblock {\em J. Pure Appl. Algebra}, 45(3):261--272, 1987.

\bibitem{Serre55}
Jean-Pierre Serre.
\newblock G\'eom\'etrie alg\'ebrique et g\'eom\'etrie analytique.
\newblock {\em Ann. Inst. Fourier, Grenoble}, 6:1--42, 1955--1956.

\bibitem{TY}
Ser~Peow Tan, Yan~Loi Wong, and Ying Zhang.
\newblock The {${\rm SL}(2,\Bbb C)$} character variety of a one-holed torus.
\newblock {\em Electron. Res. Announc. Amer. Math. Soc.}, 11:103--110
  (electronic), 2005.

\end{thebibliography}
  \bibliographystyle{plain}
  \nocite{RI90}


  \newpage

  \setcounter{section}{0}
  \renewcommand\thesection{\Alph{section}}
  \section{Appendix}

  Some formulas which appear in this article like (\ref{Formula1}) and
  (\ref{Formula2}) are really easy to prove by double induction on the
  corresponding indices, since both members of the expression satisfy
  the same recursive equation and hence the proofs are reduced to the base cases.
  However, finding those identities would require to work out all the
  computations as in the proof of Lemma \ref{demcomple}.
  In this appendix we explain how to find them and show those complicated formulas
  in detail without using induction. For the proofs given $(X,Y,Z)\in \C^3$, let us consider
  $A,B\in SL(2,\C)$ two matrices such that $(X,Y,Z) = (\tr A, \tr B, \tr AB)$.

  \begin{lemma}
  \begin{small}
  $$
    s_m(X) D =
    \begin{cases}
      \tr A^{\frac{m+1}{2}}BA^{-1}B^{-1}+\tr A^{\frac{m-1}{2}}BA^{-1}B^{-1}-\tr A^{\frac{m-3}{2}}
        -\tr A^{\frac{m-1}{2}} & \text{if $m$ is odd}, \\
      \tr A^{\frac{m}{2}}BA^{-1}B^{-1}-\tr A^{\frac{m-2}{2}} & \text{if $m$ is even}.
    \end{cases}
  $$
  \end{small}
  \end{lemma}

  \begin{proof}
  We will make use of Proposition \ref{propibas}(5), so it is clear that we must work
  out separately three different cases. Also we will assume that $m$ is a non-negative
  integer. The negative case follows immediately from the positive one together with
  the third part of Proposition \ref{propibas} and the formulas given in (\ref{proptra}).

  \vspace{1cm}
  \noindent -) $m$ odd.
  \begin{small}
  \begin{equation*}
  \begin{split}
    s_m(X) & D(X,Y,Z) = \Bigl( 1+ \sum_{i=1}^{\frac{m-1}{2}} f_i(X) \Bigr) D(X,Y,Z)
    = \Bigl( 1+ \sum_{i=1}^{\frac{m-1}{2}} \tr A^i \Bigr) \Bigl( \tr ABA^{-1}B^{-1}-2 \Bigr) =\\
    = & \ D(X,Y,Z)\ +\ \sum_{i=1}^{\frac{m-1}{2}} \Bigl[ \tr A^{i+1}BA^{-1}B^{-1}\ +\
    \tr (A^{i-1})(BAB^{-1}) \Bigr]\ -\ 2 \sum_{i=1}^{\frac{m-1}{2}} \tr A^i =\\
    = & \ D + \sum_{i=1}^{\frac{m-1}{2}} \Bigl[ \tr A^{i+1}BA^{-1}B^{-1} +
    (\tr A^{i-1} \tr A) - \tr A^{i-1}BA^{-1}B^{-1} \Bigr] - 2\sum_{i=1}^{\frac{m-1}{2}} \tr A^i =\\
    = & \ D + \sum_{i=1}^{\frac{m-1}{2}} \Bigl[\tr A^{i+1}BA^{-1}B^{-1} + (\tr A^i + \tr A^{i-2}) -
    \tr A^{i-1}BA^{-1}B^{-1}\Bigr] - 2 \sum_{i=1}^{\frac{m-1}{2}} \tr A^i = \\
    = & \ D \ +\ \sum_{i=1}^{\frac{m-1}{2}} \Bigl(\tr A^{i-2} - \tr A^i \Bigr)
    \ +\ \sum_{i=1}^{\frac{m-1}{2}} \Bigl(\tr A^{i+1}BA^{-1}B^{-1} -
    \tr A^{i-1}BA^{-1}B^{-1} \Bigr) = \\
    = & \ \Bigl(\tr ABA^{-1}B^{-1} - 2 \Bigr) + \left(\tr A^{-1} + \tr A^{0} - \tr A^{\frac{m-3}{2}} -
    \tr A^{\frac{m-1}{2}} \right) +\\
    & + \left( \tr A^{\frac{m-1}{2}}BA^{-1}B^{-1}  + \tr A^{\frac{m+1}{2}}BA^{-1}B^{-1}
    - \tr A^0BA^{-1}B^{-1} - \tr ABA^{-1}B^{-1} \right) = \\[0.1cm]
    = & \ \tr A^{\frac{m+1}{2}} BA^{-1}B^{-1} + \tr A^{\frac{m-1}{2}} BA^{-1}B^{-1}
    - \tr A^{\frac{m-3}{2}} - \tr A^{\frac{m-1}{2}}.
  \end{split}
  \end{equation*}
  \end{small}

  \noindent -) $m\equiv 0$ modulo 4.
  \begin{small}
  \begin{equation*}
  \begin{split}
    s_m(X) & D(X,Y,Z)\ =\ \Bigl(\sum_{i=1}^{\frac{m}{4}} f_{2i-1}(X)\Bigr) D(X,Y,Z) =
    \Bigl( \sum_{i=1}^{\frac{m}{4}} \tr A^{2i-1} \Bigr) \Bigl( \tr ABA^{-1}B^{-1} - 2 \Bigr) =\\
    = & \ \sum_{i=1}^{\frac{m}{4}} \Bigl[ \tr A^{2i}BA^{-1}B^{-1} + \tr (A^{2i-2})(BAB^{-1}) \Bigr]
     - 2 \sum_{i=1}^{\frac{m}{4}} \tr A^{2i-1} =\\
    = & \ \sum_{i=1}^{\frac{m}{4}} \Bigl[ \tr A^{2i}BA^{-1}B^{-1} + (\tr A^{2i-2} \tr A)
    - \tr A^{2i-2}BA^{-1}B^{-1} \Bigr] - 2 \sum_{i=1}^{\frac{m}{4}} \tr A^{2i-1} =\\
    = & \ \sum_{i=1}^{\frac{m}{4}} \Bigl[ \tr A^{2i}BA^{-1}B^{-1} + (\tr A^{2i-1} + \tr A^{2i-3})
    - \tr A^{2i-2}BA^{-1}B^{-1} \Bigr] - 2 \sum_{i=1}^{\frac{m}{4}} \tr A^{2i-1} =\\
    = & \ \sum_{i=1}^{\frac{m}{4}} \Bigl( \tr A^{2i-3} - \tr A^{2i-1} \Bigr) + \sum_{i=1}^{\frac{m}{4}}
    \Bigl( \tr A^{2i} BA^{-1}B^{-1} - \tr A^{2i-2} BA^{-1}B^{-1} \Bigr) = \\
    = & \ \Bigl(\tr A^{-1} - \tr A^{\frac{m-2}{2}}\Bigr) + \Bigl(\tr A^{\frac{m}{2}}BA^{-1}B^{-1}
    - \tr A^0BA^{-1}B^{-1} \Bigr) =\\[0.1cm]
    = & \ \tr A^{\frac{m}{2}}BA^{-1}B^{-1} - \tr A^{\frac{m-2}{2}}.
  \end{split}
  \end{equation*}
  \end{small}

  \vspace{1cm}
  \noindent -) $m\equiv 2$ modulo 4.
  \begin{small}
  \begin{equation*}
  \begin{split}
    s_m(X) & D(X,Y,Z)\ =\ \Bigl(1+\sum_{i=1}^{\frac{m-2}{4}} f_{2i}(X)\Bigr)\cdot D =
    \Bigl(1 + \sum_{i=1}^{\frac{m-2}{4}} \tr A^{2i} \Bigr)\Bigl( \tr ABA^{-1}B^{-1} -2\Bigr) =\\
    = & \ D(X,Y,Z) + \sum_{i=1}^{\frac{m-2}{4}} \Bigl[ \tr A^{2i+1}BA^{-1}B^{-1} + \tr (A^{2i-1})(BAB^{-1}) \Bigr]
    - 2 \sum_{i=1}^{\frac{m-2}{4}} \tr A^{2i} =\\
    = & \ D + \sum_{i=1}^{\frac{m-2}{4}} \Bigl[ \tr A^{2i+1}BA^{-1}B^{-1} + (\tr A^{2i-1} \tr A) -
    \tr A^{2i-1}BA^{-1}B^{-1} \Bigr] - 2 \sum_{i=1}^{\frac{m-2}{4}} \tr A^{2i} =\\
    = & \ D + \sum_{i=1}^{\frac{m-2}{4}} \Bigl[ \tr A^{2i+1}BA^{-1}B^{-1} + (\tr A^{2i} + \tr A^{2i-2}) -
    \tr A^{2i-1}BA^{-1}B^{-1} \Bigr] - 2 \sum_{i=1}^{\frac{m-2}{4}} \tr A^{2i} =\\
    = & \ D + \sum_{i=1}^{\frac{m-2}{4}} \Bigl( \tr A^{2i-2} - \tr A^{2i} \Bigr) +
    \sum_{i=1}^{\frac{m-2}{4}} \Bigl( \tr A^{2i+1}BA^{-1}B^{-1} - \tr A^{2i-1}BA^{-1}B^{-1} \Bigr) =\\
    = & \ \Bigl( \tr ABA^{-1}B^{-1} - 2 \Bigr) + \Bigl( \tr A^{0} - \tr A^{\frac{m-2}{2}} \Bigr) +
    \Bigl( \tr A^{\frac{m}{2}}BA^{-1}B^{-1} - \tr ABA^{-1}B^{-1} \Bigr) = \\[0.1cm]
    = & \ \tr A^{\frac{m}{2}}BA^{-1}B^{-1} - \tr A^{\frac{m-2}{2}}.
  \end{split}
  \end{equation*}
  \end{small}
  \end{proof}

  \begin{lemma}\label{lemmaPrevious}
  $k$ an integer.
  \begin{equation*}
  \begin{split}
  \sum_{j=1}^{\frac{n-1}{2}} & (-1)^j \tr B^j \tr A^{k}BA^{-1}B^{-1} =\\
  \stackrel{(1)}{=}\ & \sum_{j=1}^{\frac{n-1}{2}} (-1)^j \tr B^j A^{k-1} + \sum_{j=1}^{\frac{n-1}{2}}
  (-1)^j \tr B^{j+2} A^{-k+1} - \tr A^{-k+1}B + \tr A^{-k}BAB - \\
  & \ -\ (-1)^{\frac{n-3}{2}} \tr B^{\frac{n-1}{2}}A^{-k}BA \ -\
  (-1)^{\frac{n-1}{2}} \tr B^{\frac{n+1}{2}}A^{-k}BA = \\
  \stackrel{(2)}{=}\ & \sum_{j=1}^{\frac{n-1}{2}} (-1)^j \tr B^j A^{k-1} + \sum_{j=1}^{\frac{n-1}{2}}
  (-1)^j \tr B^{j-2} A^{-k+1} + \tr A^{-k+1}B^{-1} - \tr A^{-k}BAB^{-1} +\\
  & \ +\ (-1)^{\frac{n-3}{2}} \tr A^{-k}B^{\frac{n-1}{2}}AB^{-1} \ +\
  (-1)^{\frac{n-1}{2}} \tr A^{-k}B^{\frac{n+1}{2}}AB^{-1}.
  \end{split}
  \end{equation*}
  \end{lemma}

  \begin{proof}
  \begin{small}
  \begin{equation*}
  \begin{split}
  {\bf 1)}\quad \sum_{j=1}^{\frac{n-1}{2}} & (-1)^j \tr B^j \tr BA^{-1}B^{-1}A^{k} \ =\
  \sum_{j=1}^{\frac{n-1}{2}} (-1)^j \left[ \tr (B^{j+1}) (A^{-1}B^{-1}A^{k}) +
  \tr B^{j-1}A^{-k}BA \right] = \\
  =\ & \sum_{j=1}^{\frac{n-1}{2}} (-1)^j \Bigl[ \tr B^{j+1} \tr B^{-1}A^{k-1}
  - \tr B^{j+1}A^{-k}BA + \tr B^{j-1}A^{-k}BA \Bigr] =\\
  =\ & \sum_{j=1}^{\frac{n-1}{2}} (-1)^j \Bigl[ \tr B^{j} A^{k-1} + \tr B^{j+2} A^{-k+1}
  + \left(\tr B^{j-1}A^{-k}BA - \tr B^{j+1}A^{-k}BA \right) \Bigr] =\\
  =\ & \sum_{j=1}^{\frac{n-1}{2}} (-1)^j \tr B^j A^{k-1} + \sum_{j=1}^{\frac{n-1}{2}}
  (-1)^j \tr B^{j+2} A^{-k+1} - \tr B^{0}A^{-k}BA + \tr BA^{-k}BA +\\
  & \ -\ (-1)^{\frac{n-3}{2}} \tr B^{\frac{n-1}{2}}A^{-k}BA\ -\
  (-1)^{\frac{n-1}{2}} \tr B^{\frac{n+1}{2}}A^{-k}BA.\\[0.5cm]
  {\bf 2)}\quad \sum_{j=1}^{\frac{n-1}{2}} & (-1)^j \tr B^j \tr A^{k}BA^{-1}B^{-1} \ =\
  \sum_{j=1}^{\frac{n-1}{2}} (-1)^j \left[ \tr (B^{j-1}) (A^{k} BA^{-1}) +
  \tr B^{j+1}AB^{-1}A^{-k} \right] = \\
  =\ & \sum_{j=1}^{\frac{n-1}{2}} (-1)^j \Bigl[ \tr B^{j-1} \tr A^{k-1}B
  - \tr B^{j-1}AB^{-1}A^{-k}
  + \tr B^{j+1}AB^{-1}A^{-k} \Bigr] =\\
  =\ & \sum_{j=1}^{\frac{n-1}{2}} (-1)^j \Bigl[ \tr B^{j} A^{k-1} + \tr B^{j-2} A^{-k+1}
  + \left(\tr B^{j+1}AB^{-1}A^{-k} - \tr B^{j-1}AB^{-1}A^{-k}\right) \Bigr] =\\
  =\ & \sum_{j=1}^{\frac{n-1}{2}} (-1)^j \tr B^j A^{k-1} + \sum_{j=1}^{\frac{n-1}{2}}
  (-1)^j \tr B^{j-2} A^{-k+1} + \tr B^{0}AB^{-1}A^{-k} - \tr BAB^{-1}A^{-k} +\\
  & \ +\ (-1)^{\frac{n-3}{2}} \tr B^{\frac{n-1}{2}}AB^{-1}A^{-k}\ +\
  (-1)^{\frac{n-1}{2}} \tr B^{\frac{n+1}{2}}AB^{-1}A^{-k}.
  \end{split}
  \end{equation*}
  \end{small}
  \end{proof}

  \begin{lemma}
  $m$ and $n$ odd.
  \begin{equation*}
  \begin{split}
    s_m(X)\sigma_n(Y)D = & \left(F_{\frac{m+3}{2},\frac{n-1}{2}}
    - F_{\frac{m-3}{2},\frac{n+1}{2}}\right) +
    \left(F_{\frac{m+1}{2},\frac{n-1}{2}}
    - F_{\frac{m-1}{2},\frac{n+1}{2}}\right) + \\
    & + \left(F_{\frac{m-1}{2},\frac{n+3}{2}}
    - F_{\frac{m+1}{2},\frac{n-3}{2}}\right) +
    \left(F_{\frac{m-3}{2},\frac{n+3}{2}}
    - F_{\frac{m+3}{2},\frac{n-3}{2}}\right) +\\
    & + F_{1,-1} \left[\left(F_{\frac{m+1}{2},\frac{n-1}{2}}-F_{\frac{m-1}{2},
    \frac{n+1}{2}}\right) + \left(F_{\frac{m-1}{2},\frac{n-1}{2}}-F_{\frac{m+1}{2},
    \frac{n+1}{2}}\right)\right].
  \end{split}
  \end{equation*}
  \end{lemma}

  \begin{proof}
  We will use Lemma \ref{lemmaPrevious}(1) in $(*)$ taking $k=(m+1)/2$ and $k=(m-1)/2$.
  \begin{equation*}
  \begin{split}
    (-1)&^{\frac{n-1}{2}} \sigma_n(Y) s_m(X) D =
    \Bigl( 1 + \sum_{j=1}^{\frac{n-1}{2}} (-1)^j \tr B^j \Bigr)
    \left( \tr A^{\frac{m+1}{2}}BA^{-1}B^{-1}+\tr A^{\frac{m-1}{2}}BA^{-1}B^{-1} - \right. \\
    & \left. - \tr A^{\frac{m-3}{2}} - \tr A^{\frac{m-1}{2}} \right) \stackrel{(*)}{=} \\[0.25cm]
    \stackrel{(*)}{=}\ & \tr A^{\frac{m+1}{2}}BA^{-1}B^{-1}+\tr A^{\frac{m-1}{2}}BA^{-1}B^{-1}
    -\tr A^{\frac{m-3}{2}} - \tr A^{\frac{m-1}{2}} + \\
    & + \sum (-1)^j \tr B^j A^{\frac{m-3}{2}}
    + \sum (-1)^j \tr B^{j+2} A^{\frac{-m+3}{2}} - \tr A^{\frac{-m+3}{2}}B
    + \tr A^{\frac{-m+1}{2}} BAB -\\
    & - (-1)^{\frac{n-3}{2}} \tr B^{\frac{n-1}{2}} A^{\frac{-m+1}{2}} BA
    - (-1)^{\frac{n-1}{2}} \tr B^{\frac{n+1}{2}} A^{\frac{-m+1}{2}}BA +\\
    & + \sum (-1)^j \tr B^j A^{\frac{m-1}{2}}
    + \sum (-1)^j \tr B^{j+2} A^{\frac{-m+1}{2}}
    - \tr A^{\frac{-m+1}{2}}B + \tr A^{\frac{-m-1}{2}}BAB -\\
    & - (-1)^{\frac{n-3}{2}} \tr B^{\frac{n-1}{2}} A^{\frac{-m-1}{2}} BA
    - (-1)^{\frac{n-1}{2}} \tr B^{\frac{n+1}{2}} A^{\frac{-m-1}{2}}BA
    - \sum (-1)^j \tr B^j A^{\frac{m-3}{2}} - \\
    & - \sum (-1)^j \tr B^j A^{\frac{-m+3}{2}}
    - \sum (-1)^j \tr B^j A^{\frac{m-1}{2}}
    - \sum (-1)^j \tr B^j A^{\frac{-m+1}{2}} = \\[0.5cm]
    =\ & \tr A^{\frac{m+1}{2}}BA^{-1}B^{-1}+\tr A^{\frac{m-1}{2}}BA^{-1}B^{-1}
    -\tr A^{\frac{m-3}{2}} - \tr A^{\frac{m-1}{2}} - \tr A^{\frac{-m+3}{2}}B +\\
    & + \tr A^{\frac{-m+1}{2}} BAB
    - (-1)^{\frac{n-3}{2}} \tr B^{\frac{n-1}{2}} A^{\frac{-m+1}{2}} BA
    - (-1)^{\frac{n-1}{2}} \tr B^{\frac{n+1}{2}} A^{\frac{-m+1}{2}}BA -\\
    & - \tr A^{\frac{-m+1}{2}}B + \tr A^{\frac{-m-1}{2}}BAB
    - (-1)^{\frac{n-3}{2}} \tr B^{\frac{n-1}{2}} A^{\frac{-m-1}{2}} BA - \\
    & - (-1)^{\frac{n-1}{2}} \tr B^{\frac{n+1}{2}} A^{\frac{-m-1}{2}}BA
    + \sum (-1)^j \left[ \tr B^{j+2} A^{\frac{-m+3}{2}}
    - \tr B^j A^{\frac{-m+3}{2}}\right] + \\
    & + \sum (-1)^j \left[ \tr B^{j+2} A^{\frac{-m+1}{2}}
    - \tr B^j A^{\frac{-m+1}{2}}\right] = \\[0.5cm]
    =\ & \tr A^{\frac{m+1}{2}}BA^{-1}B^{-1}+\tr A^{\frac{m-1}{2}}BA^{-1}B^{-1}
    -\tr A^{\frac{m-3}{2}} - \tr A^{\frac{m-1}{2}} - \tr A^{\frac{-m+3}{2}}B +\\
    & + \tr A^{\frac{-m+1}{2}} BAB
    - (-1)^{\frac{n-3}{2}} \tr B^{\frac{n-1}{2}} A^{\frac{-m+1}{2}} BA
    - (-1)^{\frac{n-1}{2}} \tr B^{\frac{n+1}{2}} A^{\frac{-m+1}{2}}BA -\\
    & - \tr A^{\frac{-m+1}{2}}B + \tr A^{\frac{-m-1}{2}}BAB
    - (-1)^{\frac{n-3}{2}}\left[ \tr B^{\frac{n-1}{2}} A^{\frac{-m-1}{2}} BA
    - \tr B^{\frac{n+1}{2}} A^{\frac{-m-1}{2}}BA\right] + \\
    & + \left[ \tr B A^{\frac{-m+3}{2}} - \tr B^2 A^{\frac{-m+3}{2}}
    + (-1)^{\frac{n-3}{2}} \tr B^{\frac{n+1}{2}} A^{\frac{-m+3}{2}}
    + (-1)^{\frac{n-1}{2}} \tr B^{\frac{n+3}{2}} A^{\frac{-m+3}{2}} \right] +\\
    & + \left[ \tr B A^{\frac{-m+1}{2}} - \tr B^2 A^{\frac{-m+1}{2}}
    + (-1)^{\frac{n-3}{2}} \tr B^{\frac{n+1}{2}} A^{\frac{-m+1}{2}}
    + (-1)^{\frac{n-1}{2}} \tr B^{\frac{n+3}{2}} A^{\frac{-m+1}{2}} \right] =\\
  \end{split}
  \end{equation*}
  \begin{equation*}
  \begin{split}
    =\ & \tr A^{\frac{m+1}{2}}BA^{-1}B^{-1}+\tr A^{\frac{m-1}{2}}BA^{-1}B^{-1}
    -\tr A^{\frac{m-3}{2}} - \tr A^{\frac{m-1}{2}} + \tr A^{\frac{-m+1}{2}} BAB -\\
    & - (-1)^{\frac{n-3}{2}}\left[ \tr B^{\frac{n-1}{2}} A^{\frac{-m+1}{2}} \tr BA
    - \tr B^{\frac{n-3}{2}} A^{\frac{-m-1}{2}} \right]
    - (-1)^{\frac{n-1}{2}} \left[ \tr B^{\frac{n+1}{2}}A^{\frac{-m+1}{2}} \tr BA - \right.\\
    & \left. - \tr B^{\frac{n-1}{2}} A^{\frac{-m-1}{2}} \right] + \tr A^{\frac{-m-1}{2}}BAB -
    (-1)^{\frac{n-3}{2}} \left[ \tr B^{\frac{n-1}{2}} A^{\frac{-m-1}{2}} \tr BA
    - \tr B^{\frac{n-3}{2}} A^{\frac{-m-3}{2}} \right] -\\
    & - (-1)^{\frac{n-1}{2}}\left[ \tr B^{\frac{n+1}{2}} A^{\frac{-m-1}{2}} \tr BA
    - \tr B^{\frac{n-1}{2}} A^{\frac{-m-3}{2}}\right] - \tr B^2 A^{\frac{-m+3}{2}} + \\
    & + (-1)^{\frac{n-3}{2}} \tr B^{\frac{n+1}{2}} A^{\frac{-m+3}{2}}
    + (-1)^{\frac{n-1}{2}} \tr B^{\frac{n+3}{2}} A^{\frac{-m+3}{2}} - \tr B^2 A^{\frac{-m+1}{2}} + \\
    & + (-1)^{\frac{n-3}{2}} \tr B^{\frac{n+1}{2}} A^{\frac{-m+1}{2}}
    + (-1)^{\frac{n-1}{2}} \tr B^{\frac{n+3}{2}} A^{\frac{-m+1}{2}} = \\[0.5cm]
    =\ & (-1)^{\frac{n-1}{2}} \tr AB \left[ \left( \tr A^{\frac{-m-1}{2}} B^{\frac{n-1}{2}}
    - \tr A^{\frac{-m+1}{2}} B^{\frac{n+1}{2}}\right) +
    \left( \tr A^{\frac{-m+1}{2}} B^{\frac{n-1}{2}}
    - \tr A^{\frac{-m-1}{2}} B^{\frac{n+1}{2}}\right) \right] + \\
    & + (-1)^{\frac{n-1}{2}} \left[ \left( \tr A^{\frac{-m-3}{2}} B^{\frac{n-1}{2}}
    - \tr A^{\frac{-m+3}{2}} B^{\frac{n+1}{2}}\right) +
    \left( \tr A^{\frac{-m-1}{2}} B^{\frac{n-1}{2}}
    - \tr A^{\frac{-m+1}{2}} B^{\frac{n+1}{2}}\right) \right] + \\
    & + (-1)^{\frac{n-1}{2}} \left[ \left( \tr A^{\frac{-m+1}{2}} B^{\frac{n+3}{2}}
    - \tr A^{\frac{-m-1}{2}} B^{\frac{n-3}{2}}\right) +
    \left( \tr A^{\frac{-m+3}{2}} B^{\frac{n+3}{2}}
    - \tr A^{\frac{-m-3}{2}} B^{\frac{n-3}{2}}\right) \right] + \\
    & + \left( \tr A^{\frac{m+1}{2}}BA^{-1}B^{-1} + \tr A^{\frac{-m-1}{2}}BAB
    - \tr A^{\frac{m-1}{2}} - \tr A^{\frac{-m+1}{2}}B^2 \right) +\\
    & + \left( \tr A^{\frac{m-1}{2}}BA^{-1}B^{-1} + \tr A^{\frac{-m+1}{2}} BAB
    -\tr A^{\frac{m-3}{2}} - \tr A^{\frac{-m+3}{2}}B^2 \right) \stackrel{(**)}{=} \\[0.5cm]
    =\ & (-1)^{\frac{n-1}{2}} \tr AB \left[ \left( \tr A^{\frac{-m-1}{2}} B^{\frac{n-1}{2}}
    - \tr A^{\frac{-m+1}{2}} B^{\frac{n+1}{2}}\right) +
    \left( \tr A^{\frac{-m+1}{2}} B^{\frac{n-1}{2}}
    - \tr A^{\frac{-m-1}{2}} B^{\frac{n+1}{2}}\right) \right] + \\
    & + (-1)^{\frac{n-1}{2}} \left[ \left( \tr A^{\frac{-m-3}{2}} B^{\frac{n-1}{2}}
    - \tr A^{\frac{-m+3}{2}} B^{\frac{n+1}{2}}\right) +
    \left( \tr A^{\frac{-m-1}{2}} B^{\frac{n-1}{2}}
    - \tr A^{\frac{-m+1}{2}} B^{\frac{n+1}{2}}\right) \right] + \\
    & + (-1)^{\frac{n-1}{2}} \left[ \left( \tr A^{\frac{-m+1}{2}} B^{\frac{n+3}{2}}
    - \tr A^{\frac{-m-1}{2}} B^{\frac{n-3}{2}}\right) +
    \left( \tr A^{\frac{-m+3}{2}} B^{\frac{n+3}{2}}
    - \tr A^{\frac{-m-3}{2}} B^{\frac{n-3}{2}}\right) \right] = \\[0.5cm]
    =\ & (-1)^{\frac{n-1}{2}} F_{1,-1} \left[ \left(F_{\frac{m+1}{2},\frac{n-1}{2}}-F_{\frac{m-1}{2},
    \frac{n+1}{2}}\right) + \left(F_{\frac{m-1}{2},\frac{n-1}{2}}-F_{\frac{m+1}{2},
    \frac{n+1}{2}}\right)\right] + \\
    & + (-1)^{\frac{n-1}{2}} \left[ \left(F_{\frac{m+3}{2},\frac{n-1}{2}}
    - F_{\frac{m-3}{2},\frac{n+1}{2}}\right) +
    \left(F_{\frac{m+1}{2},\frac{n-1}{2}}
    - F_{\frac{m-1}{2},\frac{n+1}{2}}\right)\right] + \\
    & + (-1)^{\frac{n-1}{2}} \left[ \left(F_{\frac{m-1}{2},\frac{n+3}{2}}
    - F_{\frac{m+1}{2},\frac{n-3}{2}}\right) +
    \left(F_{\frac{m-3}{2},\frac{n+3}{2}}
    - F_{\frac{m+3}{2},\frac{n-3}{2}}\right) \right].\\
  \end{split}
  \end{equation*}

  \vspace{0.1cm}

  \begin{equation*}
  \hspace{1cm} (**)\hspace{-1cm} \begin{split}
  \tr A^k B A^{-1}B^{-1} =\ & \tr B A^{-1}B^{-1}A^{k} =
  \tr B \tr A^{k-1}B^{-1} - \tr BA^{-k}BA =\\
  =\ & \tr A^{k-1} + \tr A^{-k+1}B^2 - \tr A^{-k}BAB.
  \end{split}
  \end{equation*}

  We have used $(**)$ taking $k=(m+1)/2$ and $k=(m-1)/2$.
  \end{proof}

\newpage

  \begin{lemma}
  $m$ even and $n$ odd.
  \begin{equation*}
  \begin{split}
    s_m(X)\sigma_n(Y)D = & \left(F_{\frac{m+2}{2},\frac{n+1}{2}}
    - F_{\frac{m-2}{2},\frac{n-1}{2}}\right) +
    \left(F_{\frac{m-2}{2},\frac{n-3}{2}}
    - F_{\frac{m+2}{2},\frac{n+3}{2}}\right) + \\
    & + F_{1,1} \left(F_{\frac{m}{2},\frac{n+1}{2}}-F_{\frac{m}{2},
    \frac{n-1}{2}}\right).
  \end{split}
  \end{equation*}
  \end{lemma}

  \begin{proof}
  We will use Lemma \ref{lemmaPrevious}(2) in $(*)$ taking $k=m/2$.
  \begin{equation*}
  \begin{split}
    (-1)^{\frac{n-1}{2}} & \sigma_n(Y) s_m(X) D = \Bigl( 1 + \sum_{j=1}^{\frac{n-1}{2}} (-1)^j \tr B^j \Bigr)
    \left( \tr A^{\frac{m}{2}}BA^{-1}B^{-1} - \tr A^{\frac{m-2}{2}} \right) =\\
    =\ & \tr A^{\frac{m}{2}}BA^{-1}B^{-1} - \tr A^{\frac{m-2}{2}} + \sum_{j=1}^{\frac{n-1}{2}} (-1)^j
    \left[ \tr B^j \tr A^{\frac{m}{2}}BA^{-1}B^{-1} - \tr B^j \tr A^{\frac{m-2}{2}}\right]\stackrel{(*)}{=}\\[0.25cm]
    \stackrel{(*)}{=}\ & \tr A^{\frac{m}{2}}BA^{-1}B^{-1} - \tr A^{\frac{m-2}{2}} + \sum
    (-1)^j \tr B^j A^{\frac{m-2}{2}} + \sum (-1)^j \tr B^{j-2} A^{\frac{-m+2}{2}} +\\
    & + \tr A^{\frac{-m+2}{2}}B^{-1} - \tr A^{\frac{-m}{2}}BAB^{-1} + (-1)^{\frac{n-3}{2}}
    \tr A^{\frac{-m}{2}} B^{\frac{n-1}{2}} AB^{-1} +\\
    & + (-1)^{\frac{n-1}{2}} \tr A^{\frac{-m}{2}}B^{\frac{n+1}{2}}AB^{-1}
    - \sum (-1)^j \tr B^j A^{\frac{m-2}{2}} -
    \sum (-1)^j \tr B^jA^{\frac{-m+2}{2}} =\\[0.5cm]
    =\ & \tr A^{\frac{m}{2}}BA^{-1}B^{-1} - \tr A^{\frac{m-2}{2}}
    + \tr A^{\frac{-m+2}{2}}B^{-1} - \tr A^{\frac{-m}{2}}BAB^{-1} +\\
    & + (-1)^{\frac{n-3}{2}} \tr A^{\frac{-m}{2}} B^{\frac{n-1}{2}} AB^{-1} +
    (-1)^{\frac{n-1}{2}} \tr A^{\frac{-m}{2}}B^{\frac{n+1}{2}}AB^{-1} - \tr B^{-1}A^{\frac{-m+2}{2}} + \\
    & + \tr B^{0} A^{\frac{-m+2}{2}} - (-1)^{\frac{n-3}{2}} \tr B^{\frac{n-3}{2}} A^{\frac{-m+2}{2}}
    - (-1)^{\frac{n-1}{2}} \tr B^{\frac{n-1}{2}} A^{\frac{-m+2}{2}} =\\[0.5cm]
    =\ & (-1)^{\frac{n-3}{2}} \tr A^{\frac{-m}{2}}B^{\frac{n-1}{2}} \tr AB^{-1}
    - (-1)^{\frac{n-3}{2}} \tr A^{\frac{-m-2}{2}}B^{\frac{n+1}{2}} +\\
    & + (-1)^{\frac{n-1}{2}} \tr A^{\frac{-m}{2}}B^{\frac{n+1}{2}} \tr AB^{-1}
    - (-1)^{\frac{n-1}{2}} \tr A^{\frac{-m-2}{2}}B^{\frac{n+3}{2}} -\\
    & - (-1)^{\frac{n-3}{2}} \tr B^{\frac{n-3}{2}} A^{\frac{-m+2}{2}}
    - (-1)^{\frac{n-1}{2}} \tr B^{\frac{n-1}{2}} A^{\frac{-m+2}{2}} =\\[0.5cm]
    =\ & (-1)^{\frac{n-1}{2}} \left[ \left( \tr A^{\frac{-m-2}{2}} B^{\frac{n+1}{2}} -
    \tr A^{\frac{-m+2}{2}} B^{\frac{n-1}{2}} \right)
    + \left( \tr A^{\frac{-m+2}{2}} B^{\frac{n-3}{2}} -
    \tr A^{\frac{-m-2}{2}} B^{\frac{n+3}{2}} \right) + \right.\\
    & \left. \ +\ \tr AB^{-1} \left( \tr A^{\frac{-m}{2}} B^{\frac{n+1}{2}} -
    \tr A^{\frac{-m}{2}} B^{\frac{n-1}{2}} \right) \right] =\\[0.5cm]
    =\ & (-1)^{\frac{n-1}{2}} \left[ \left(F_{\frac{m+2}{2},\frac{n+1}{2}}
    - F_{\frac{m-2}{2},\frac{n-1}{2}}\right) +
    \left(F_{\frac{m-2}{2},\frac{n-3}{2}}
    - F_{\frac{m+2}{2},\frac{n+3}{2}}\right) +\right. \\
    & \left. \ +\ F_{1,1} \left(F_{\frac{m}{2},\frac{n+1}{2}}-F_{\frac{m}{2},
    \frac{n-1}{2}}\right) \right].
  \end{split}
  \end{equation*}
  \end{proof}

\newpage

  \begin{lemma}
  $m$ odd, $n$ even.
  \begin{equation*}
  \begin{split}
    s_m(X)f_{\frac{n}{2}}(Y)D=&\left(F_{\frac{m+3}{2},\frac{n}{2}}-F_{\frac{m-3}{2},
    \frac{n}{2}}\right) + \left(F_{\frac{m+1}{2},\frac{n}{2}}-F_{\frac{m-1}{2},\frac{n}{2}}\right) + \\
    &+ \left(F_{\frac{m-1}{2},\frac{n-4}{2}}-F_{\frac{m+1}{2},
    \frac{n+4}{2}}\right) +
    \left(F_{\frac{m-3}{2},\frac{n-4}{2}}-F_{\frac{m+3}{2},\frac{n+4}{2}}\right) + \\
    &+ F_{1,1} \left[\left(F_{\frac{m+1}{2},\frac{n+2}{2}}-F_{\frac{m-1}{2},
    \frac{n-2}{2}}\right)+\left(F_{\frac{m-1}{2},\frac{n+2}{2}}-F_{\frac{m+1}{2},
    \frac{n-2}{2}}\right)\right].
  \end{split}
  \end{equation*}
  \end{lemma}

  \begin{proof}
    \begin{equation*}
      \begin{split}
        f_{\frac{n}{2}}(Y)s_m&(X)D =\\[0.25cm] =\ & \tr B^{\frac{n}{2}}\left(\tr A^{\frac{m+1}{2}}BA^{-1}B^{-1}+
       \tr A^{\frac{m-1}{2}}BA^{-1}B^{-1}-\tr A^{\frac{m-3}{2}}-\tr A^{\frac{m-1}{2}}\right) = \\[0.5cm]
=\ &\tr B^{\frac{n-2}{2}}A^{\frac{m+1}{2}}BA^{-1}+\tr A^{\frac{-m-1}{2}}B^{\frac{n-2}{2}}AB^{-1}+\tr B^{\frac{n-2}{2}}A^{\frac{m-1}{2}}BA^{-1}+\\
&+\tr A^{\frac{-m+1}{2}}B^{\frac{n+2}{2}}AB^{-1}-\tr B^{\frac{n}{2}}A^{\frac{m-3}{2}}-\tr B^{\frac{n}{2}}A^{\frac{-m+3}{2}}-\tr B^{\frac{n}{2}}A^{\frac{m-1}{2}}-\\&
-\tr B^{\frac{n}{2}}A^{\frac{-m+1}{2}}=\\[0.5cm]
=\ &\tr B^{\frac{n}{2}}\tr A^{\frac{m-1}{2}}B-\tr A^{\frac{-m-1}{2}}B^{\frac{n-2}{2}}AB^{-1}+\tr A^{\frac{-m-1}{2}}B^{\frac{n+2}{2}}AB^{-1}+\\&+\tr B^{\frac{n}{2}}\tr A^{\frac{m-3}{2}}B-\tr A^{\frac{-m+1}{2}}B^{\frac{n-2}{2}}AB^{-1}+\tr A^{\frac{-m+1}{2}}B^{\frac{n+2}{2}}AB^{-1}-\\&-\tr B^{\frac{n}{2}}A^{\frac{m-3}{2}}-\tr B^{\frac{n}{2}}A^{\frac{-m+3}{2}}-\tr B^{\frac{n}{2}}A^{\frac{m-1}{2}}-\tr B^{\frac{n}{2}}A^{\frac{-m+1}{2}}=\\[0.5cm]
=\ &\tr B^{\frac{n-4}{2}}A^{\frac{-m+1}{2}}-\tr A^{\frac{-m-1}{2}}B^{\frac{n-2}{2}}\tr AB^{-1}+\tr A^{\frac{-m-3}{2}}B^{\frac{n}{2}}+\\
&+\tr A^{\frac{-m-1}{2}}B^{\frac{n+2}{2}}\tr AB^{-1}-\tr A^{\frac{-m-3}{2}}B^{\frac{n+4}{2}}+\tr B^{\frac{n-4}{2}}A^{\frac{-m+3}{2}}-\\
&-\tr A^{\frac{-m+1}{2}}B^{\frac{n-2}{2}}\tr AB^{-1}+\tr A^{\frac{-m-1}{2}}B^{\frac{n}{2}}+\tr A^{\frac{-m+1}{2}}B^{\frac{n+2}{2}}\tr AB^{-1}-\\&-\tr A^{\frac{-m-1}{2}}B^{\frac{n+4}{2}}-\tr B^{\frac{n}{2}}A^{\frac{-m+3}{2}}-\tr B^{\frac{n}{2}}A^{\frac{-m+1}{2}}=\\[0.5cm]
=\ &\left(\tr A^{\frac{-m-3}{2}}B^{\frac{n}{2}}-\tr A^{\frac{-m+3}{2}}B^{\frac{n}{2}}\right)+\left(\tr A^{\frac{-m-1}{2}}B^{\frac{n}{2}}-\tr A^{\frac{-m+1}{2}}B^{\frac{n}{2}}\right)+\\
&+\left(\tr A^{\frac{-m+1}{2}}B^{\frac{n-4}{2}}-\tr A^{\frac{-m-1}{2}}B^{\frac{n+4}{2}}\right)+\left(\tr A^{\frac{-m+3}{2}}B^{\frac{n-4}{2}}-\tr A^{\frac{-m-3}{2}}B^{\frac{n+4}{2}}\right)+\\&+\tr AB^{-1}\left[\left(\tr A^{\frac{-m-1}{2}}B^{\frac{n+2}{2}}-\tr A^{\frac{-m+1}{2}}B^{\frac{n-2}{2}}\right)\right.+\\&+\left.\left(\tr A^{\frac{-m+1}{2}}B^{\frac{n+2}{2}}-\tr A^{\frac{-m-1}{2}}B^{\frac{n-2}{2}}\right)\right]=\\[0.5cm]
        =\ &\left(F_{\frac{m+3}{2},\frac{n}{2}}-F_{\frac{m-3}{2},
        \frac{n}{2}}\right) + \left(F_{\frac{m+1}{2},\frac{n}{2}}-F_{\frac{m-1}{2},\frac{n}{2}}\right) + \\
        &+ \left(F_{\frac{m-1}{2},\frac{n-4}{2}}-F_{\frac{m+1}{2},
        \frac{n+4}{2}}\right) +
        \left(F_{\frac{m-3}{2},\frac{n-4}{2}}-F_{\frac{m+3}{2},\frac{n+4}{2}}\right) + \\
        &+ F_{1,1} \left[\left(F_{\frac{m+1}{2},\frac{n+2}{2}}-F_{\frac{m-1}{2},
        \frac{n-2}{2}}\right)+\left(F_{\frac{m-1}{2},\frac{n+2}{2}}-F_{\frac{m+1}{2},
        \frac{n-2}{2}}\right)\right].
      \end{split}
    \end{equation*}
  \end{proof}

  \begin{lemma}
  $m$, $n$ even.
  \begin{equation*}
  \begin{split}
    s_m(X)f_{\frac{n}{2}}(Y)D=&\left(F_{\frac{m+2}{2},\frac{n}{2}}-F_{\frac{m-2}{2},
    \frac{n}{2}}\right) + \left(F_{\frac{m-2}{2},\frac{n-4}{2}}-F_{\frac{m+2}{2},\frac{n+4}{2}}\right) + \\
    & + F_{1,1} \left(F_{\frac{m}{2},\frac{n+2}{2}}-F_{\frac{m}{2},
    \frac{n-2}{2}}\right).
  \end{split}
  \end{equation*}
  \end{lemma}

  \begin{proof}
  \begin{equation*}
  \begin{split}
    f_{\frac{n}{2}}(Y) & s_m(X)D \ =\
    \tr B^{\frac{n}{2}} \left( \tr A^{\frac{m}{2}}BA^{-1}B^{-1} - \tr A^{\frac{m-2}{2}} \right) = \\[0.25cm]
    =\ & \tr B^{\frac{n-2}{2}} A^{\frac{m}{2}}BA^{-1} + \tr B^{\frac{n+2}{2}}AB^{-1}A^{\frac{-m}{2}}
    - \tr B^{\frac{n}{2}} A^{\frac{m-2}{2}} - \tr B^{\frac{n}{2}} A^{\frac{-m+2}{2}} = \\[0.25cm]
    =\ & \tr B^{\frac{n-2}{2}} \tr A^{\frac{m-2}{2}} B - \tr A^{\frac{-m}{2}} B^{\frac{n-2}{2}}AB^{-1}
    + \tr A^{\frac{-m}{2}}B^{\frac{n+2}{2}}AB^{-1} - \\
    & - \tr B^{\frac{n}{2}} A^{\frac{m-2}{2}} - \tr B^{\frac{n}{2}} A^{\frac{-m+2}{2}} = \\[0.25cm]
    =\ & \tr B^{\frac{n}{2}} A^{\frac{m-2}{2}} + \tr B^{\frac{n-4}{2}} A^{\frac{-m+2}{2}}
    - \tr A^{\frac{-m}{2}} B^{\frac{n-2}{2}} \tr AB^{-1} + \tr A^{\frac{-m-2}{2}} B^{\frac{n}{2}} + \\
    & + \tr A^{\frac{-m}{2}} B^{\frac{n+2}{2}} \tr AB^{-1} - \tr A^{\frac{-m-2}{2}} B^{\frac{n+4}{2}}
    - \tr B^{\frac{n}{2}} A^{\frac{m-2}{2}} - \tr B^{\frac{n}{2}} A^{\frac{-m+2}{2}} = \\[0.25cm]
    =\ & \left(\tr A^{\frac{-m-2}{2}} B^{\frac{n}{2}} - \tr A^{\frac{-m+2}{2}} B^{\frac{n}{2}}\right)
    + \left( \tr A^{\frac{-m+2}{2}}B^{\frac{n-4}{2}} - \tr A^{\frac{-m-2}{2}} B^{\frac{n+4}{2}} \right) +\\
    & + \tr AB^{-1} \left( \tr A^{\frac{-m}{2}}B^{\frac{n+2}{2}}
    -\tr A^{\frac{-m}{2}}B^{\frac{n-2}{2}}\right)=\\[0.25cm]
    =\ & \left(F_{\frac{m+2}{2},\frac{n}{2}}-F_{\frac{m-2}{2},\frac{n}{2}}\right)
    + \left(F_{\frac{m-2}{2},\frac{n-4}{2}}-F_{\frac{m+2}{2},\frac{n+4}{2}}\right)
    + F_{1,1} \left(F_{\frac{m}{2},\frac{n+2}{2}}-F_{\frac{m}{2},\frac{n-2}{2}}\right).
  \end{split}
  \end{equation*}
  \end{proof}

  %

\end{document}